\documentclass[reqno]{amsart}
\usepackage[utf8]{inputenc}
\usepackage{pgffor}
\usepackage{xcolor,soul}
\definecolor{hyperlink}{rgb}{0.7 0 0}
\usepackage[colorlinks=true,allcolors=hyperlink]{hyperref}
\usepackage{booktabs}
\usepackage{mathtools}
\usepackage{mathdots}
\usepackage{amsthm}
\usepackage{amssymb}
\usepackage{graphicx}
\usepackage{cases}
\usepackage{array}
\usepackage{enumitem}
\usepackage[figureposition=top, skip=3pt]{caption}
\usepackage{subcaption}
\usepackage{placeins}
\usepackage[style=authoryear,
            sorting=nyt,
            backend=biber,
            natbib,
            minbibnames=7,
            maxbibnames=10]{biblatex}
\usepackage{tabularx, array}
\usepackage{booktabs,array}
\usepackage{array}
\usepackage{multirow}

\graphicspath{ {./images/} }

\usepackage[vlined,boxed,commentsnumbered]{algorithm2e} 
\usepackage{listings}
\usepackage[scaled]{beramono}
\usepackage[T1]{fontenc}
\definecolor{codegreen}{rgb}{0,0.6,0}
\definecolor{codegrey}{rgb}{0.5,0.5,0.5}
\definecolor{codepurple}{rgb}{0.58,0,0.82}
\definecolor{backcolor}{rgb}{0.95,0.95,0.92}
\newcommand{\code}[1]{\texttt{\fontsize{9.5pt}{10pt}\selectfont #1}}
\newcommand{\smallcode}[1]{\texttt{\fontsize{7.5pt}{8pt}\selectfont #1}}

\usepackage{tikz}
\usetikzlibrary{decorations.pathmorphing, patterns,shapes}
\newcommand{\zigzag}[1]{\begin{tikzpicture}
  \draw[decoration = {zigzag, segment length = 3mm, amplitude = 0.5mm},decorate, color = codegrey] (0,0.1)--(#1,0.1);
\end{tikzpicture}
}

\usepackage{mdframed}

\newmdenv[
  topline=false,
  bottomline=false,
  rightline=false,
  skipabove=\topsep,
  skipbelow=\topsep
]{leftrule}

\usepackage{siunitx}
\newcolumntype{N}{S[scientific-notation = true, round-mode = figures, round-precision = 3, table-format=2.2e-1]}

\newcolumntype{K}{S[
    table-format           = 6.0,        group-separator        = {,},        table-number-alignment = center,     table-text-alignment   = center    ]}

\newcolumntype{T}[1]{  S[
    table-format       = #1.0 ,
    group-separator    = {,} ,
    group-minimum-digits = 4
  ]}

\usepackage{array,siunitx}
\newcolumntype{J}[2]{S[
  table-format = #1.#2,
  round-mode   = places,
  round-precision = #2
]}

\makeatletter

\DeclareCaptionSubType*{figure}
\theoremstyle{definition}

\theoremstyle{remark}

\renewcommand{\vec}{\mathbf}
\renewcommand*\d{\mathop{}\!\mathrm{d}}

\newcommand{\customlabel}[2]{#2\def\@currentlabel{#2}\label{#1}}

\newcommand{\varparallel}{\mathbin{\!/\mkern-5mu/\!}}

\DeclarePairedDelimiter{\Norm}{\Vert}{\Vert}
\DeclarePairedDelimiter{\norm}{|}{|}
\DeclarePairedDelimiter{\inner}{\langle}{\rangle}

\DeclareCiteCommand{\citeyear}
    {}
    {\bibhyperref{\printdate}}
    {\multicitedelim}
    {}

\DeclareCiteCommand{\citeyearpar}
    {}
    {\mkbibparens{\bibhyperref{\printdate}}}
    {\multicitedelim}
    {}
\makeatother

\calclayout

\begin{document}
\raggedbottom
\title[Fitting an escalier to a curve]
    {    Fitting an escalier to a curve}

\begin{abstract}
    We analyze the problem of fitting a \emph{fonction en escalier} or multi-step function to a curve in $L^2$ Hilbert space.  We propose a two-stage optimization approach whereby the step positions are initially fixed, corresponding to a classic linear least-squares problem with closed-form solution, and then are allowed to vary, leading to first-order conditions that can be solved recursively. We find that, subject to regularity conditions, the speed of convergence is linear as the number of steps $n$ goes to infinity, and we develop a simple algorithm to recover the global optimum fit.  Our numerical results based on a sweep search implementation show promising performance in terms of speed and accuracy.
\end{abstract}

\author{S\'ebastien Bossu\textsuperscript{*}} 
\thanks{*~UNC Charlotte Department of Mathematics and Statistics. SB is the corresponding author and thanks Stephan Sturm for useful discussions in early versions of this paper, and Joseph Bakita for his earlier involvement with the numerical section.  NE joined the collaboration in 2025 and has been in charge of the numerical section~\ref{sec:numerics}. The research of SB and NE is partially supported by funds provided by The University of North Carolina at Charlotte. }
\author{Andrew Papanicolaou\textsuperscript{\dag}}
\thanks{\dag~NCSU Department of Mathematics. Andrew Papanicolaou is partially supported by NSF grant DMS-1907518
\\
\indent Our C++ code is available online at \url{https://github.com/sbossu/EscalierSweep}}
\author{Nour El Hatto\textsuperscript{*}}

\date{\today}

\maketitle



\section{Introduction}  \label{sec:intro}

\subsection{Problem statement}

In the Hilbert space $L^2([a,b])$ of square-integrable functions over a finite interval $[a,b]\subset\mathbb R$ with canonical inner product $\inner{f,g} = \int_a^b f(x)g(x) \d x$ and induced norm $\Norm{f} = \sqrt{\inner{f,f}}$, given a target function $f\in L^2([a,b])$, find a \emph{fonction en escalier} $g$ (``stairway function'', also known as \emph{step function}, \emph{simple function}, or \emph{piecewise continuous function}) with $n\geq 1$ steps that minimizes the $L^2$ distance
\begin{equation}    \label{eq:L2-min-problem}
    \min_{g \in E_n} \Norm{ f - g }
    =
    \min_{\boldsymbol\phi,\vec k} \Norm*{ f - \phi_0 - \sum_{i=1}^n \phi_i\, u_{k_i} },   
\end{equation}
where $ u_{k} $ denotes the unit step function centered at point $k\in[a,b]$, and $ E_n \subset L^2([a,b]) $ is the subset of $n$-step escalier functions of the form
\[
    g(x) = \phi_0 + \sum_{i=1}^n \phi_i\, u_{k_i}(x),
    \qquad x\in[a,b],
\]
for some parameters $ (\phi_i)_{0\leq i\leq n} \in \mathbb R^{n+1}, (k_i)_{1\leq i\leq n} \in \mathbb (a,b)^n$.  Without loss of generality we assume $ a =: k_0 < k_1 < \cdots < k_n < k_{n+1} \coloneqq b$ throughout this paper, and we define the vectors of $n+1$ coefficients with zero-based indexing $ \boldsymbol\phi = (\phi_i)_{0\leq i\leq n}, \vec k = (k_i)_{0\leq i\leq n}$.
Because the vector $\boldsymbol \phi$ is unconstrained in $\mathbb R^{n+1}$, the least-squares optimization problem \eqref{eq:L2-min-problem} may be split into a linear subproblem followed by a non-linear super-problem as
\begin{equation}    \label{eq:L2-min-problem-sub-super}
    \min_{\boldsymbol\phi,\vec k} \Norm*{ f - \sum_{i=0}^n \phi_i\, u_{k_i} } = \min_{\vec k} \min_{\boldsymbol\phi} \Norm*{f - \sum_{i=0}^n \phi_i u_{k_i} }.
\end{equation}

\subsection{Motivation and organization of the paper}

The problem of finding an optimal escalier fit to a given curve is relevant to many areas of applied mathematics, including:
\begin{itemize}[leftmargin=*,parsep=6pt]
    \item Machine learning: escalier fitting corresponds to training a single-layer neural network with step activation functions \citep[e.g.][]{szandala:2021}.  While sigmoid and ReLU activation functions are more popular choices, step functions are very easy to compute and correspond to the boolean architecture of computers.  They may be a suitable choice for binary classification algorithms. By the \emph{universal approximation theorem}, a wide class of target functions $f(x)$ may be approximated by a one-hidden-layer network with non-polynomial activation function to an arbitrary degree of precision---see \citet{guilhoto2018overview} for a mathematical overview of neural networks.  The step function used in this paper is the simplest kind of activation function, which is the limit-case of a sigmoid and the derivative of the ubiquitous Rectified Linear Unit (ReLU) activation function.
    \item Wavelet analysis: escalier fitting is an alternative to a Haar basis decomposition of the target function $f$.  The Haar basis is orthogonal but nonparametric and numerous Haar wavelets may be required to achieve the desired precision of approximation. For general background on wavelets, see \citet{mallat:2009}, \citet{donoho:1994}, \citet{candes:1999:curvelets}.
    \item Spline method and Finite Element Method (FEM): escalier fitting may be viewed as a special case of splining and FEM.  FEM has become a popular method for numerically solving differential equations arising in engineering and mathematical modeling.  For a mathematical treatment of FEM, see \citet{strang-fix:2008}.
    \item Quadrature methods: escalier fitting corresponds to an optimal partition of the interval $[a,b]$ for Riemann sum approximation of $\int_a^b f$,  which bears some resemblance with Gaussian quadrature whose knots $k_1 < \cdots < k_n$ are optimal with respect to a class of orthogonal polynomials.
\end{itemize}
 Further background on nonlinear approximation using orthonormal bases in $L^p$ Hilbert spaces can be found in the survey by \citet{devore:1998}.  While our approach is limited to the $L^2$ case for ease of exposition and derivation, we use a nonorthonormal basis which is more general in nature and better suited to machine learning applications, and more importantly we propose an algorithm to compute the global optimum fit.

The rest of this paper is organized as follows: Section~2 derives the optimal escalier fit when the step positions $k_1 < \cdots < k_n$ are fixed and investigates its asymptotic properties.  Section~3 derives first-order conditions for the step positions when they are allowed to vary and the target function $f$ is differentiable, and finds that speed of convergence is $O(1/n)$ subject to general conditions.  Section~4, which is the main contribution of this paper, presents our algorithm to compute the global optimum fit based on recurrence relations stemming from the first-order conditions derived in Section~3.  Section~5 presents our C++ implementation and testing results for a selection of target functions with various characteristics.

\section{Linear least-squares escalier fit}

In this section, we assume that the step positions $ a =: k_0 < k_1 < \cdots < k_n < k_{n+1} \coloneqq b$ are fixed parameters, and we review the solution to the linear least-squares subproblem in equation \eqref{eq:L2-min-problem-sub-super},
\[
    \min_{\boldsymbol \phi} \Norm*{f - \sum_{i=0}^n \phi_i u_{k_i} }.
\]

\subsection{Optimal fit}

\begin{figure}[htbp]
    \centering
    \includegraphics[width=0.75\linewidth]{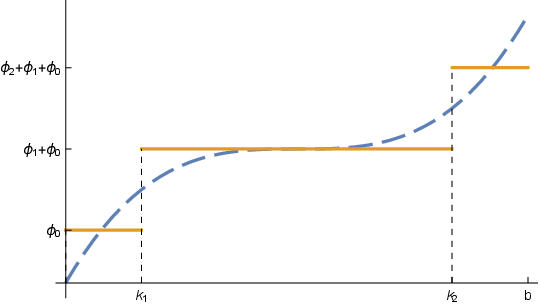}
    \caption{Two-step optimal escalier fit to a curve}
    \label{fig:escalier-2step-example}
\end{figure}

 By standard linear algebra theory \citep[e.g.][ch.10, pp.171--219]{muscat:2014}, the least-squares solution is unique and given as
\begin{equation}  \label{eq:lsq-sol}
    \phi_i^* = \sum_{j=0}^n g^{i,j} \inner{f,u_{k_j}},
    \qquad 0\leq i\leq n,
\end{equation}
where $(g^{i,j})_{0\leq i,j\leq n} =: \vec G^{-1} $ is the inverse of the $(n+1)\times(n+1)$ Gram matrix $\vec G \coloneqq \left( \inner{u_{k_i},u_{k_j}} \right)_{0\leq i,j\leq n}$ generated by the system of step functions $(u_{k_0},\ldots,u_{k_n})$.  Alternatively, in vector notations,
\begin{equation}  \label{eq:lsq-sol-vec}
    \boldsymbol{\phi}^* = \vec G^{-1} \vec f,
    \qquad \vec f \coloneqq \left(\inner{f,u_{k_0}}, \cdots, \inner{f,u_{k_n}}\right)^T.
\end{equation}
The Gram coefficients are given as
\[
    g_{i,j} \coloneqq \inner{u_{k_i}, u_{k_j}} = \int_a^b u_{k_i}(x)u_{k_j}(x) \d x
    = b - k_{\max(i,j)},
    \qquad 0 \leq i,j \leq n,
\]
so that the Gram matrix is a \emph{single-pair matrix} (see \cite{gantmacher:2002}, \cite{baranger-ducjacquet}, \cite{meurant:1992}, \cite{vandebril:2008}) with tridiagonal inverse \citep[see][]{bossu:2024}:
\begin{equation}    \label{eq:BC-Gram-inv}
    \vec G^{-1} = \begin{pmatrix}
        \frac{1}{k_1-a} & \frac{-1}{k_1-a} &  &  & (0) \\
     \frac{-1}{k_1- a} & \frac{k_2-a}{\left(k_1-a\right) \left(k_2-k_1\right)} & \ddots & &  \\
      & \ddots & \ddots & \ddots &  \\
      & & \ddots & \frac{k_{n}-k_{n-2}}{\left(k_{n-1} - k_{n-2}\right) \left(k_n - k_{n-1}\right)} & \frac{-1}{k_n-k_{n-1}} \\
      (0) & & & \frac{-1}{k_n-k_{n-1}} & \frac{(b-k_{n-1})/(b-k_n)}{k_n - k_{n-1}}
    \end{pmatrix}.
\end{equation}
Note that each core diagonal coefficient is the opposite sum of its surrounding row coefficients:
\[
    g^{i,i} = \frac{k_{i+1}-k_{i-1}}{(k_{i+1}-k_{i})(k_{i}-k_{i-1})}
    = \frac1{k_{i+1}-k_{i}} + \frac1{k_{i}-k_{i-1}}
    = -g^{i,i+1} - g^{i,i-1},
    \qquad 1 \leq i \leq n -1.
\]
After simplification, optimal quantities $\boldsymbol \phi^* = \vec G^{-1} \vec f$ are given as
\begin{equation}        \label{eq:binary-lsq-sol}
    \begin{dcases}
        \phi_0^* = \frac{\inner{f,u_a-u_{k_1}}}{k_1 - a}
        \\
        \phi_i^* = \frac{\inner{f,u_{k_i}-u_{k_{i+1}}}}{k_{i+1} - k_{i}} - \frac{\inner{f,u_{k_{i-1}}-u_{k_i}}}{k_i - k_{i-1}},
        & 1 \leq i \leq n.
    \end{dcases}
\end{equation}
Recognizing
\[
    \frac{\inner{f,u_{k_{i-1}}-u_{k_i}}}{k_i - k_{i-1}}
    =
    \frac1{k_i - k_{i-1}} \int_a^b f(x)\big[u_{k_{i-1}}(x)-u_{k_{i}}(x)\big]\! \d x
    = 
    \frac1{k_i - k_{i-1}} \int_{k_{i-1}}^{k_i} f(x) \d x
\]
as the mean value of the target function $f$ over the interval $[k_{i-1}, k_{i}]$, we see that the optimal fit is an escalier where the height of each step is the mean value of $f$ over each interval, as illustrated in Figure \ref{fig:escalier-2step-example}.

\subsection{Quality of fit}

At the optimum, the squared norm of the residual is known to be
\begin{equation}        \label{eq:TSS-ESS}
    \Norm*{f - \sum_{i=0}^n \phi_i^* u_{k_i} }^2 = \Norm f^2 - \Norm*{\sum_{i=0}^n \phi_i^* u_{k_i} }^2,    
\end{equation}
which is the linear algebra abstract expression of the statistical equation SSE = TSS -- ESS.  Expanding the squares; switching to vector notation with the help of the Gram matrix $\vec G \coloneqq \left( \inner{u_{k_i},u_{k_j}} \right)_{0\leq i,j\leq n}$; then substituting equation \eqref{eq:lsq-sol-vec} and simplifying,
\[
    \text{ESS} = 
    \Norm*{\sum_{i=0}^n \phi_i^* u_{k_i} }^2
    = \sum_{0\leq i,j\leq n} \phi_i^*\phi_j^* \inner*{ u_{k_i}, u_{k_j}}
    = \boldsymbol \phi^{*T} \vec G \boldsymbol \phi^*
    = \vec{f}^T \vec G^{-1} \vec{f},
\]
so that the coefficient of determination (R-squared) of the function approximation is simply
\begin{equation}        \label{eq:lsq-r2}
    \text R^2 = \frac{\text{ESS}}{\text{TSS}} = \frac {\vec f^T \vec G^{-1} \vec f}{\Norm f^2}.    
\end{equation}
Substituting the tridiagonal coefficients of equation \eqref{eq:BC-Gram-inv},
\begin{align*}
    \text{ESS} = \vec{f}^T \vec G^{-1} \vec{f}
    & =
    \sum_{0\leq i,j\leq n} g^{i,j} \inner{f,u_{k_i}}\inner{f,u_{k_j}}
    \\
    & =
    \frac{\inner{f, u_a}^2}{k_1 - a} + 
    \sum_{i=1}^n \frac{\inner{f, u_{k_i}}^2 (k_{i+1}-k_{i-1})}{\left(k_i - k_{i-1}\right) \left(k_{i+1}-k_i\right)}
    - 2 \sum_{i=1}^{n} \frac{\inner{f, u_{k_{i-1}}}\inner{f,u_{k_i}}}{k_i-k_{i-1}}.
\end{align*}
Substituting $\frac{k_{i+1}-k_{i-1}}{\left(k_i - k_{i-1}\right) \left(k_{i+1}-k_i\right)} = \frac1{k_i - k_{i-1}} + \frac1{k_{i+1}-k_i}$, reindexing one sum and recognizing a perfect square,
\begin{equation}        \label{eq:binary-ESS}
    \text{ESS} = \vec{f}^T \vec G^{-1} \vec{f}
    =
    \sum_{i=1}^{n+1} \frac{\inner{f,u_{k_i} -  u_{k_{i-1}}}^2}{k_i-k_{i-1}}
    .
\end{equation}
Substituting $ \inner{f,u_{k_i} -  u_{k_{i-1}}} = -\int_{k_{i-1}}^{k_i} f(x) \d x $ and simplifying,
\begin{align}
    \text{ESS}
    & = \vec{f}^T \vec G^{-1} \vec{f}
    = \sum_{i=1}^{n+1} \frac1{k_i-k_{i-1}}\left(\int_{k_{i-1}}^{k_i} f(x) \d x \right)^2
    \label{eq:binary-ESS3}
    \\
    & = \sum_{i=1}^{n+1} (k_i-k_{i-1})\!\left(\frac{\int_{k_{i-1}}^{k_i} f(x) \d x}{k_i-k_{i-1}}\right)^2.
    \label{eq:binary-ESS4}
\end{align}
Observe that equation \eqref{eq:binary-ESS4} is a weighted average of squared mean values of $f$ which is numerically more stable than equation \eqref{eq:binary-ESS3}.

\subsection{Asymptotics}

Equation \eqref{eq:binary-ESS4} corresponds to the Riemann-type sum
\[
    \vec{f}^T \vec G^{-1} \vec{f}
    =
    \sum_{i=1}^{n+1} \inner*{f,\frac{u_{k_i} - u_{k_{i-1}}}{k_i-k_{i-1}}}^2(k_i-k_{i-1})
\]
Intuitively, as $n\to\infty$ and the width $k_i-k_{i-1} \to 0$, we have $\frac{u_{k_i} - u_{k_{i-1}}}{k_i-k_{i-1}} \longrightarrow \frac{\partial u_{k_{i-1}}}{\partial k_{i-1}} =\delta_{k_{i-1}}$ which is Dirac's delta function centered at $ x = k_{i-1}$, and we recover
\[
    \vec{f}^T \vec G^{-1} \vec{f}
    \longrightarrow \int_a^b f(x)^2 \d x \eqqcolon \Norm{f}^2
\]
as expected.  In terms of optimal quantities $\boldsymbol\phi^*$ from equation \eqref{eq:binary-lsq-sol}, we have
\begin{equation*}       
    \begin{dcases}
        \phi_0^* = \inner*{f,\frac{u_a-u_{k_1}}{k_1 - a}} \longrightarrow \inner{f, \delta_a} = f(a)
        \\
        \phi_i^* = \frac{\inner{f,u_{k_i}-u_{k_{i+1}}}}{k_{i+1} - k_{i}} - \frac{\inner{f,u_{k_{i-1}}-u_{k_i}}}{k_i - k_{i-1}} \longrightarrow 0,
        & 1 \leq i \leq n \to \infty.
    \end{dcases}
\end{equation*}
Assuming $k_{i+1}-k_i = k_i - k_{i-1} \eqqcolon \Delta k > 0 $, we have the second-order central finite difference.
\[
    \frac{\frac{u_{k_i}-u_{k_{i+1}}}{k_{i+1} - k_{i}} - \frac{u_{k_{i-1}}-u_{k_i}}{k_i - k_{i-1}}}{\Delta k} = \frac{ 2 u_{k_i}-u_{k_i+\Delta k} - u_{k_i - \Delta k} }{\Delta k^2} \longrightarrow -\frac{\partial^2 u_{k_i}}{\partial k_i^2} = - \delta'_{k_i}.
\]
Asymptotically,
\[
    \phi_i^* \to \inner{f,-\delta'_{k_i}} \d k_i = \inner{f',\delta_{k_i}} \d k_i = f'(k_i) \d k_i,
    \qquad 1 \leq i \leq n \to \infty,
\]
provided that $f$ is differentiable. Thence
\[
    \sum_{i=0}^n \phi_i^* u_{k_i} \approx f(a) + \sum_{i=1}^n (k_{i+1} - k_i) f'(k_i) u_{k_i} 
    \longrightarrow f(a) + \int_a^b f'(k) u_k \d k = \int_a^b f(k)\delta_k \d k = f
\]
as expected, where we integrated by parts and used the sifting property of Dirac's delta function in the last steps.

\subsection{Error bound}        \label{sec:error-bound}

Substituting \eqref{eq:binary-ESS3}, the error may be rewritten as
\begin{align*}
    \text{SSE} = \text{TSS}-\text{ESS}
    & = \int_a^b f(x)^2 \d x - \sum_{i=1}^{n+1} \frac1{k_i-k_{i-1}}\left(\int_{k_{i-1}}^{k_i} f(x) \d x \right)^2
    \\
    & = \sum_{i=1}^{n+1} \int_{k_{i-1}}^{k_i} f(x) \frac{\int_{k_{i-1}}^{k_i} [f(x) - f(y)] \d y}{k_i-k_{i-1}} \d x.
\end{align*}
Assuming $f$ is differentiable and $f'$ is bounded, by the mean-value inequality we have $ \norm{f(x) - f(y)} \leq \norm{x-y} \sup_{(k_{i-1}, k_i)} | f' | $.  By triangular inequality,
\[
    \text{SSE}
    \leq \sum_{i=1}^{n+1} \sup_{(k_{i-1}, k_i)} \norm{ f' } \cdot \int_{k_{i-1}}^{k_i} \norm{ f(x) } \frac{\int_{k_{i-1}}^{k_i} \norm{x - y} \d y}{k_i-k_{i-1}} \d x.
\]
It is easy to show that $ \frac{\int_{k_{i-1}}^{k_i} \norm{x - y} \d y}{k_i-k_{i-1}} \leq \frac{k_i - k_{i-1}}2 $, whence
\[
    \text{SSE}
    \leq \sum_{i=1}^{n+1}  \frac{k_i - k_{i-1}}2 \sup_{(k_{i-1}, k_i)} \norm{ f' } \cdot \int_{k_{i-1}}^{k_i} \norm{ f(x) } \d x
    \leq \frac12 \sup_{(a,b)} \norm{ f } \cdot \sup_{(a,b)} \norm{ f' } \cdot \sum_{i=1}^{n+1} (k_i - k_{i-1})^2.
\]
Thus, if $ k_i - k_{i-1} = O(1/n)$ as $n \to \infty$ then $ \text{SSE} = O(1/n) $.

\section{Nonlinear least-squares escalier fit}      \label{sec:nonlinear-lsq}

Having solved the linear regression component of problem \eqref{eq:L2-min-problem}, in this section we allow the step positions $ k_1 < \dots < k_n$ to vary between $a =: k_0$ and $b =: k_{n+1}$, and we emphasize dependence on $\vec k \coloneqq (k_0, k_1, \dots, k_n)$ with a subscript or functional notation. We study the nonlinear super-problem
\begin{equation}        \label{eq:L2-min-problem-nlsq}
    \min_{\vec k} \Norm*{f - \sum_{i=0}^n \phi_i^*(\vec k) u_{k_i} }
\end{equation}
where $ (\phi_i^*(\vec k))_{0\leq i\leq n} =: \boldsymbol \phi^*(\vec k) $ is the solution to the linear least-squares subproblem given by equation \eqref{eq:lsq-sol-vec},
\begin{equation}        \label{eq:lsq-sol-vec2}
    \boldsymbol \phi^*(\vec k) \coloneqq \vec G^{-1}_{\vec k}\vec f_{\vec k},
    \qquad
    \vec f_{\vec k} \coloneqq \left(\inner{f,u_{k_0}}, \cdots, \inner{f,u_{k_n}}\right)^T.
\end{equation}
Equation \eqref{eq:L2-min-problem-nlsq} is a nontrivial $n$-dimensional minimization problem which is nonlinear and non-convex.

\subsection{First-order conditions}

By the partition identity \eqref{eq:TSS-ESS}, the  minimization super-problem \eqref{eq:L2-min-problem-nlsq} is equivalent to the ESS (or R-squared) maximization problem 
\[
    \max_{\vec k}\, \vec f_{\vec k}^T \vec G^{-1}_{\vec k} \vec{f_k}.
\]
Alternatively, in sigma notation,
\begin{equation}        \label{eq:nlsq-sol}
    \max_{\vec k} \sum_{1\leq i,j\leq n} g^{i,j}_{\vec k}\: \inner{f,u_{k_i}}\: \inner{f,u_{k_j}},
\end{equation}
where $g^{i,j}_{\vec k}$ are the coefficients of the inverse Gram matrix \eqref{eq:BC-Gram-inv}.  Based on equations \eqref{eq:binary-ESS3}--\eqref{eq:binary-ESS4}, the gradient exists if and only if the mean values are differentiable, i.e. $f$ must be continuous. First-order conditions corresponding to the maximization problem \eqref{eq:nlsq-sol} yield after   some algebra the system of critical point equations
\begin{equation}    \label{eq:critical-points-1}
    \left[\frac{\int_{k_i}^{k_{i+1}} f(x) \d x}{k_{i+1}-k_i} - f(k_i) \right]^2 = \left[\frac{\int_{k_{i-1}}^{k_i} f(x) \d x}{k_i-k_{i-1}} - f(k_i) \right]^2,
        \qquad 1 \leq i \leq n.
\end{equation}
Equivalently,
\[
    \frac{\int_{k_i}^{k_{i+1}} f(x) \d x}{k_{i+1}-k_i} = \frac{\int_{k_{i-1}}^{k_i} f(x) \d x}{k_i-k_{i-1}}
    \quad\text{or}\quad
    \frac{\frac{\int_{k_i}^{k_{i+1}} f(x) \d x}{k_{i+1}-k_i} + \frac{\int_{k_{i-1}}^{k_i} f(x) \d x}{k_i-k_{i-1}}}2 = f(k_i),
    \qquad 1 \leq i \leq n.
\]
We may safely ignore the degenerate case $ \frac{\int_{k_i}^{k_{i+1}} f(x) \d x}{k_{i+1}-k_i} = \frac{\int_{k_{i-1}}^{k_i} f(x) \d x}{k_i-k_{i-1}} $ which implies that the step function $ u_{k_i} $ has zero coefficient $\phi_{i}^* = \frac{\int_{k_i}^{k_{i+1}} f(x) \d x}{k_{i+1}-k_i} - \frac{\int_{k_{i-1}}^{k_i} f(x) \d x}{k_i-k_{i-1}} = 0$, i.e. we have two successive steps of equal height.  Hence, nondegenerate critical points $(k_1, \ldots, k_n)$ satisfy the system of $n$ conditions in $n$ unknowns
\begin{equation}        \label{eq:critical-points-2}
    \frac{\int_{k_i}^{k_{i+1}} f(x) \d x}{k_{i+1}-k_i} \neq \frac{\int_{k_{i-1}}^{k_i} f(x) \d x}{k_i-k_{i-1}}
    \quad\text{and}\quad
    \frac{\frac{\int_{k_i}^{k_{i+1}} f(x) \d x}{k_{i+1}-k_i} + \frac{\int_{k_{i-1}}^{k_i} f(x) \d x}{k_i-k_{i-1}}}2 = f(k_i),
    \qquad 1 \leq i \leq n,
\end{equation}
which admits at least one solution since $f\in L^2([a,b])$ has bounded TSS.  \S\ref{sec:quadratic-fitting} below gives an example where the system can be solved in closed form, however in most cases it must be solved numerically.  In Section \ref{sec:escalier-algos} we propose a basic numerical algorithm to do so.

\subsection{Example: fitting a two-step escalier to a parabola}   \label{sec:quadratic-fitting}

We derive the best two-step escalier fit to the parabola $f(x) = x^2$ over the interval $[0,b]$.  By equation \eqref{eq:critical-points-2} for $n=2$ and $a = 0$, we must solve the system of 2 equations with 2 unknowns $0 < k_1 < k_2 < b$
\begin{align*}
    & \left\{\begin{aligned}
    \frac{\int_{k_1}^{k_2} x^2 \d x}{k_2-k_1} & = 2 k_2^2 - \frac{\int_{k_2}^{b} x^2 \d x}{b-k_2} \\
    \frac{\int_{0}^{k_1} x^2 \d x}{k_1} & = 2 k_1^2 - \frac{\int_{k_1}^{k_2} x^2 \d x}{k_2-k_1} = 2 k_1^2 - 2k_2^2 + \frac{\int_{k_2}^{b} x^2 \d x}{b-k_2},
    \quad 0 < k_1 < k_2 < b,
    \end{aligned}    \right. \nonumber \\
\end{align*}
\begin{subnumcases}{\text{i.e. }}
    \frac13(k_1^2 + k_1 k_2 + k_2^2) = 2 k_2^2 - \frac13(k_2^2 + k_2 b + b^2)
    \label{eq:parabola-criteq-1}
    \\
    \frac13 k_1^2 = 2 k_1^2 - 2 k_2^2 + \frac13(k_2^2 + k_2 b + b^2),
    \quad 0 < k_1 < k_2 < b,
    \label{eq:parabola-criteq-2}
\end{subnumcases}
Simplifying and solving equation \eqref{eq:parabola-criteq-1} for $k_1\in(0,k_2)$ we find that the only feasible solution is 
\[
    k_1 = \frac12\left(\sqrt{17 k_2^2-4 b k_2 - 4 b^2}-k_2\right),
    \qquad \frac{1+\sqrt{17}}8 b < k_2 < b.
\]
Substituting into equation \eqref{eq:parabola-criteq-2} and simplifying,
\[
    35 k_2^2 - k_2 \left(5 \sqrt{-4 b^2-4 b k_2+17 k_2^2}+8 b\right) - 8 b^2 = 0,
    \qquad \frac{1+\sqrt{17}}8 b < k_2 < b,
\]
whose only feasible solution is $ k_2 = \frac{b}{160} \left(\sqrt{17}+\sqrt{366 \sqrt{17}+7906}+23\right) \approx 0.7760\: b$.  Hence the optimal escalier fit is given by
\[
    \begin{cases}
        k_1 = \frac{b}{160} \left(5 + 3 \sqrt{17}+\sqrt{350 \sqrt{17}+2418}\right) \approx 0.4969\: b
        \\
        k_2 = \frac{b}{160} \left(23 + \sqrt{17}+\sqrt{366 \sqrt{17}+7906}\right) \approx 0.7760\: b.
    \end{cases}
\]

\subsection{Speed of convergence}

In general, the critical point equation \eqref{eq:critical-points-2} may be written as
\[
    \frac{\int_{k_i}^{k_{i+1}} f(x) \d x}{k_{i+1}-k_i} -f(k_i) = f(k_i) - \frac{\int_{k_{i-1}}^{k_i} f(x)\d x}{k_i-k_{i-1}} ,
        \qquad 1 \leq i \leq n.
\]
Assuming $f$ is differentiable, a second-order Taylor expansion on both sides yields as $k_{i+1}-k_i \to 0$ and $k_i - k_{i-1} \to 0$,
\[
    \frac12 f'(k_i)(k_{i+1}-k_i) + o(k_{i+1}-k_i) = \frac12 f'(k_i)(k_i-k_{i-1}) + o(k_i-k_{i-1}) ,
        \qquad 1 \leq i \leq n.
\]
If $f'(k_i)\neq 0$, we obtain that $k_{i+1}-k_i \sim k_i - k_{i-1}$ for all $1 \leq i \leq n$. The case $f'(k_i) = 0$ can be circumvented with a higher-order Taylor expansion (if $f$ is further differentiable).  Consequently, $k_{i+1}-k_i \sim \frac{b-a}n $ as $n\to\infty$.  By \S\ref{sec:error-bound} this implies that the speed of SSE convergence is $O(1/n)$.

\section{Pseudo-code algorithm to find the best escalier fit}
\label{sec:escalier-algos}

\subsection{Recurrence relations}

The first-order conditions \eqref{eq:critical-points-2} define a three-term implicit recurrence: given initial terms $ k_0 \coloneqq a $ and $ a < k_1 < b $, each successive term $ k_2,\ldots, k_n $ may be determined by solving the forward\footnote{The equation can be reversed to obtain a backward recurrence initialized at $ k_{n+1} \coloneqq b, a < k_n < b$.} recurrence relation
\[
     \frac{\int_{k_i}^{k_{i+1}} f(x) \d x}{k_{i+1}-k_i} = 2 f(k_i) - \frac{\int_{k_{i-1}}^{k_i} f(x) \d x}{k_i-k_{i-1}},
     \qquad k_i < k_{i+1} < b,\quad 1 \leq i \leq n - 1,
\]
If there is no solution $ k_{i+1} \in (k_i, b) $ then we do not have a valid critical point.  It is worth emphasizing that the equation may have multiple solutions, in which case the recurrence must be branched.  Finally, the last points $ k_{n-1}, k_n $ should also satisfy the terminal condition for $ i = n $,
\[
    \frac{\int_{k_{n-1}}^{k_n} f(x) \d x}{k_n-k_{n-1}} + \frac{\int_{k_n}^{b} f(x) \d x}{b-k_n} = 2 f(k_n).
\]
The best escalier fit then solves the one-dimensional maximization problem
\[
    \max_{a < k_1 < b} \max_{(k_2, \ldots, k_n)\in \mathcal T(k_1)} \sum_{i=1}^{n+1} (k_i - k_{i-1}) \left(\frac{\int_{k_{i-1}}^{k_i} f(x) \d x}{k_i-k_{i-1}}\right)^2
\]
where $ \mathcal T(k_1) = \{ (k_2, \ldots, k_n) : a < k_1 < \ldots < k_n < b $ \text{ satisfy equation \eqref{eq:critical-points-2}} \} is the set of all valid tails of critical points that start with $ k_1 $.  Note that if $f$ is $K$-Lipschitz over $(a,b)$, then $ k_{i+1} > k_i + \frac1K\norm*{2 f(k_i) - \frac{\int_{k_{i-1}}^{k_i} f(x) \d x}{k_i-k_{i-1}}}$, which can be used to reduce the search interval for $k_{i+1}$.

\subsection{Algorithm}

Algorithm~\ref{alg:escalier} presents the pseudo-code for maximizing the ESS of an $n$-step escalier $(k_1, \dots, k_n)$ by letting the root $ k_1 $ vary inside the interval $(a,b)$ while the tail $(k_2, \dots, k_n)$ is calculated by recurrence.  The best escalier fit is then given by the root $k_1^*$ and corresponding tail $(k_2^*, \dots, k_n^*)$ with highest ESS.

\newcommand{\maxess}{\mathit{maxess}}
\newcommand{\ess}{\mathit{ess}}
\newcommand{\maxsteps}{\mathit{maxsteps}}
\begin{algorithm}[!htbp]
\caption{Escalier pseudo-code}
\label{alg:escalier}
\small
\DontPrintSemicolon
\SetAlgoLined
\SetKwInOut{AlgName}{name}
\SetKwInOut{Input}{input}\SetKwInOut{Output}{output}
\AlgName{Escalier}
\Input{Target function $f$, interval $[a,b]$, maximum number of steps $n$}
\Output{Optimal escalier fit $(k_1^*, k_2^*, \dots, k_n^*)$ } 
\;
\KwSty{let} $\maxess = 0$
\;
\ForEach{$k_1 \in (a,b)$}{
    \KwSty{let} $\ess = \FuncSty{buildESS}(a, k_1, n)$
    \;
    \If{$\ess>\maxess$}{
        \KwSty{let} $\maxess = \ess, k_1^* = k_1$
    }
}
\KwRet $\FuncSty{getTail}(a, k_1^*, n)$
\;
\;
\SetKwProg{func}{function}{}{end}
\func{\FuncSty{buildESS}$(k_{i-1},k_i,\maxsteps)$}{
    \KwSty{let} $ \maxess = (b - k_i)\Big(\frac{\int_{k_{i}}^{b} f(x) \d x}{b-k_{i}}\Big)^2 \quad $
    \tcp{ESS of one-step tail $(k_i,b)$}
    
    \If{$ \maxsteps > 1$}{
        Solve for $ k_{i+1} \in (k_i,b) $:
        $ \displaystyle
             \frac{\int_{k_i}^{k_{i+1}} f(x) \d x}{k_{i+1}-k_i} = 2 f(k_i) - \frac{\int_{k_{i-1}}^{k_i} f(x) \d x}{k_i-k_{i-1}}
        $

        \ForEach{solution $k_{i+1}$}{
            \KwSty{let} $\maxess = \max(\maxess, \FuncSty{buildESS}(k_i,k_{i+1},\maxsteps-1))$        
                                }
    }
                            \KwRet $ \maxess + (k_i - k_{i-1})\Big(\frac{\int_{k_{i-1}}^{k_i} f(x) \d x}{k_i-k_{i-1}}\Big)^2 $
}
\func{\FuncSty{getTail}($k_{i-1},k_i,\maxsteps)$}{
    \eIf{$ \maxsteps = 1$}{
        \KwRet $k_i$
    }{
        Solve for $ k_{i+1} \in (k_i,b) $:
        $ \displaystyle
             \frac{\int_{k_i}^{k_{i+1}} f(x) \d x}{k_{i+1}-k_i} = 2 f(k_i) - \frac{\int_{k_{i-1}}^{k_i} f(x) \d x}{k_i-k_{i-1}}
        $
        \;
        \KwSty{if} \emph{no solution} \KwSty{then return} empty list \O
        \;
        \KwSty{let} $\maxess = 0,\ t^*=\, $\O
        \;
        \ForEach{solution $k_{i+1}$}{
            \KwSty{let} $t = \FuncSty{getTail}(k_i,k_{i+1},\maxsteps-1)$
            \;
            Calculate ESS of $t$
            \;
            \If{ESS of t $\,> \maxess$}{
                \KwSty{let} $\maxess = $ ESS of t,$\ t^* = t$
            }
        }
        \KwRet the concatenated list $(k_i, t^*)$
    }
}
\end{algorithm}

\subsection{Other recurrence relations}

By iteration we have the implicit full recurrence relation
\begin{equation}        \label{eq:full-rec}
     \frac{\int_{k_i}^{k_{i+1}} f(x) \d x}{k_{i+1}-k_i} 
     = 2\sum_{j=1}^i (-1)^{i-j} f(k_j) + (-1)^i \frac{\int_{a}^{k_1} f(x) \d x}{k_1-a},
     \qquad k_i < k_{i+1} < b,\quad 1 \leq i \leq n - 1,
\end{equation}
and a valid critical point must satisfy the terminal condition
\begin{equation}        \label{eq:full-rec-terminal}
     \frac{\int_{k_n}^{b} f(x) \d x}{b-k_n} = 2\sum_{j=1}^n (-1)^{n-j} f(k_j) + (-1)^n \frac{\int_{a}^{k_1} f(x) \d x}{k_1-a}.
\end{equation}
At the expense of a more convoluted expression, we may eliminate dependence of equation \eqref{eq:full-rec} on $ \int_{a}^{k_1} f(x) \d x $ which varies with $k_1$ by means of the conspicuous identity
\[
    \sum_{i=0}^{n} (k_{i+1} - k_{i}) \frac{\int_{k_i}^{k_{i+1}} f(x) \d x}{k_{i+1}-k_i}  = \int_a^b f(x) \d x.
\]
Substituting equations \eqref{eq:full-rec} and \eqref{eq:full-rec-terminal} into the above,
\begin{align*}
    \int_a^b f(x) \d x
    & =
    \int_{a}^{k_1} f(x) \d x + \sum_{i=1}^{n} (k_{i+1} - k_{i}) \left( 2\sum_{j=1}^i (-1)^{i-j} f(k_j) + (-1)^i \frac{\int_{a}^{k_1} f(x) \d x}{k_1-a}\right)
    \\
    & =
    \frac{ \int_{a}^{k_1} f(x) \d x}{k_1 - a} \sum_{i=0}^n (-1)^i (k_{i+1} - k_i)
    + 2 \sum_{j=1}^n (-1)^j f(k_j) \sum_{i=j}^n (-1)^i (k_{i+1} - k_i),
\end{align*}
where we collected terms and swapped sums in the last step.  Denoting $ S_j(\vec k) \coloneqq \sum_{i=j}^n (-1)^i (k_{i+1} - k_i)$ for any $ 0 \leq j \leq n $,
\[
    \frac{ \int_{a}^{k_1} f(x) \d x}{k_1 - a} = \frac{ \int_a^b f(x) \d x - 2 \sum_{j=1}^n (-1)^j f(k_j) S_j(\vec k)}{S_0(\vec k)},
    \qquad S_0(\vec k) \neq 0.
\]
Substituting the above into \eqref{eq:full-rec} and rearranging,
\begin{flalign}        \label{eq:full-rec-v2}
     \frac{\int_{k_i}^{k_{i+1}} f(x) \d x}{k_{i+1}-k_i} 
     =\ \mathrlap{ \frac{(-1)^i}{S_0(\vec k)} \int_{a}^{b} f(x) \d x + 2\sum_{j=1}^i (-1)^{i-j} f(k_j)\left(1 - \frac{S_j(\vec k)}{S_0(\vec k)}\right) - 2\sum_{j=i+1}^n (-1)^{i-j} f(k_j)\frac{S_j(\vec k)}{S_0(\vec k)} ,}
     \nonumber
     \\
     && S_0(\vec k) \neq 0, 0 \leq i \leq n.
\end{flalign}
where the border cases $ i = 0, i = n $ are included for completeness.

\section{Numerical implementation: sweeping the escalier}
\label{sec:numerics}

We implemented Algorithm \ref{alg:escalier} using a sweep search of the solutions $k_{i+1}$, together with a recursion. Our implementation requires that the computation of both $f(x)$ and the mean value $\bar f(a,b) \coloneqq \frac{\int_a^b f(x) \d x}{b-a} $ is available.  The main source code is provided in Appendix, while our full source code can be found online at \url{https://github.com/sbossu/EscalierSweep}.

We evaluated the performance of our algorithm, Escalier, in two categories: speed, by measuring wall clock runtime in milliseconds; and quality of solution, by measuring accuracy using the coefficient of determination $R^2$. Performance was assessed across a selection of eight functions with various characteristics (see Table~\ref{tab:testfunctions}). In low dimension ($n=2,3$), we benchmarked Escalier against a baseline exhaustive search, hereafter called Brute Force. In higher dimension ($n=2,\dots,10$), we analyzed Escalier's scalability in terms of speed and accuracy as $n$ increased. While wall clock runtime may be difficult to reproduce since it is tied to a specific combination of hardware and software, our interest lies not in absolute runtimes but rather in relative comparisons: Escalier against Brute Force in low dimension, and runtime growth in higher dimension.

\begin{table}[htbp]
\centering
\caption{Test function zoo}
\label{tab:testfunctions}
\begin{tabular}{c >{\raggedright\arraybackslash}p{5cm} >{\raggedright\arraybackslash}m{3cm} >{\centering\arraybackslash}m{5cm}}
\toprule
\textbf{Group} & \textbf{Function} & \textbf{Comment} & \textbf{Graph} \\
\midrule
\addlinespace[5pt]
II &
\( f_1(x) =  x^2 \) & 
    Smooth \& monotonic (see \S\ref{sec:quadratic-fitting}) & 
    \includegraphics[width=0.19\textwidth]{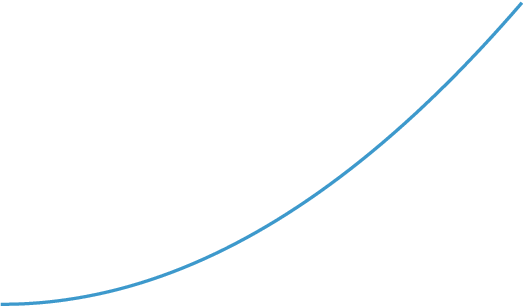} \\ \midrule 
\addlinespace[5pt]  

\addlinespace[5pt]  
I &
\( f_2(x) = \frac{1 }{\sqrt{|1-x|}} \) & 
    Non-square-integrable with midpoint singularity & 
    \includegraphics[width=0.19\textwidth]{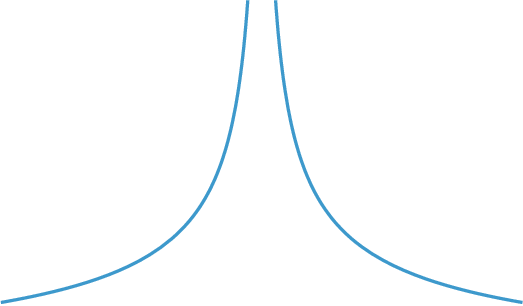} \\ \midrule 
\addlinespace[5pt]  

\addlinespace[5pt]  
II &
\( f_3(x) = \log(x)\) & 
    Monotonic, square-integrable with endpoint singularity &
    \includegraphics[width=0.19\textwidth]{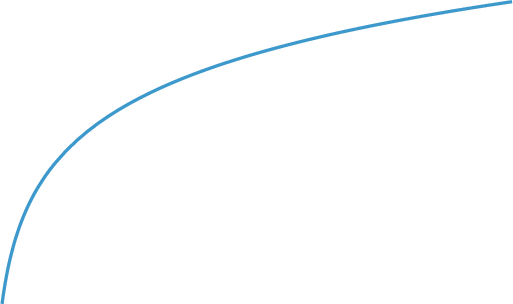} \\ \midrule 
\addlinespace[5pt]  

\addlinespace[5pt]      
I &
\( f_4(x) =  
    \begin{cases}
        0, & \text{if } x < 1 \\
        1, & \text{if } x \geq 1 
    \end{cases} \) &
    Unit step (sanity check) & 
    \includegraphics[width=0.19\textwidth]{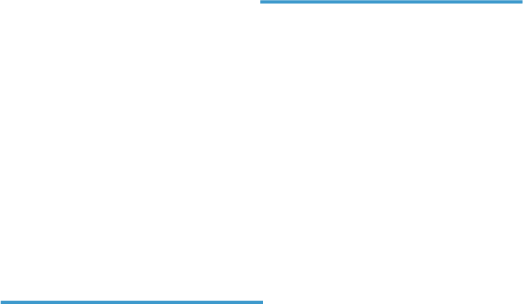} \\ \midrule 
\addlinespace[5pt]  

\addlinespace[5pt]      
I &
\( f_5(x) = \frac{1}{\sqrt[3]{x-1}}\) &
    Square-integrable with midpoint singularity &
    \includegraphics[width=0.19\textwidth]{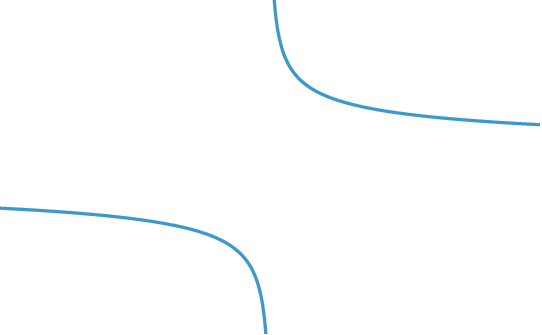} \\ \midrule 
\addlinespace[5pt]  

\addlinespace[5pt]      
II &
\( f_6(x) = \frac{1}{32} (11-10x)^2 \, (3+\cos(10\pi x))+\sin^2\!\left(\frac{\pi+10\pi x}{4}\right)\) & 
    1-D L\'evy function, smooth and highly oscillatory &
    \includegraphics[width=0.19\textwidth]{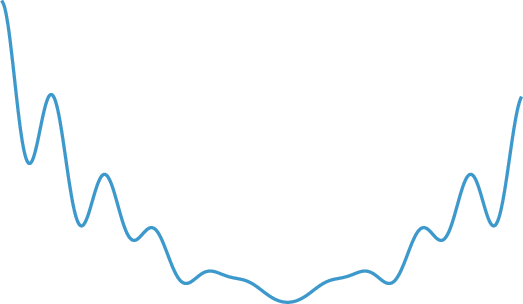} \\ \midrule 
\addlinespace[5pt]  

\addlinespace[5pt]
II &
\( f_7(x) = x\sin{(\frac{9}{4} x^2)} \) & 
    Smooth \& oscillatory &
    \includegraphics[width=0.19\textwidth]{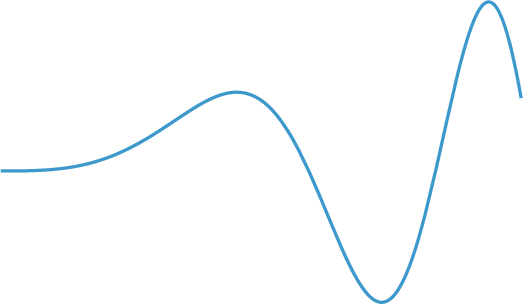} \\ \midrule

\addlinespace[5pt]
I &
\( f_8(x) =  
    \begin{cases}
        x, & \text{if } x < 1 \\
        x-1, & \text{if } x \geq 1 
    \end{cases} \)& 
    Piecewise linear with a jump discontinuity &
    \includegraphics[width=0.19\textwidth]{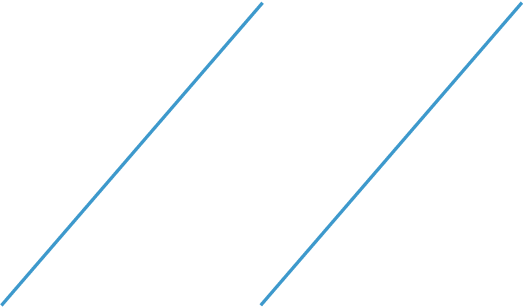} \\

\bottomrule
\end{tabular}
\end{table}

All testing was carried out on a workstation equipped with an Intel(R) Xeon(R) W5-2445 processor (3.10 GHz) and 64 GB of RAM. The code was compiled with GCC 14.2.0. Details about the compilation flags can be found in the project repository. Within this environment, experiments were conducted on the fixed interval \([a,b]=[0,2]\). Escalier was evaluated across tolerance levels \( 10^{-4}, 10^{-3}, 10^{-2}, 10^{-1}, 1\), precision levels \( 10^{-5},10^{-4}, 10^{-3}, 10^{-2}\), and maximum step count ($N_{\max}$) of 2 to 10.  Brute Force, which does not use a tolerance parameter, was run with the same precision levels \(10^{-4}, 10^{-3}, 10^{-2}\) in low dimension \(n=2,3\) only, since its exponential complexity rendered its computational cost impractical for $n > 3$.

\subsection{Implementation details}
Our C++ sweep search implementation of the pseudocode relies on some specific choices. In this subsection, we walk through the implementation block by block, highlighting the structure and any key differences.  Throughout the remainder of this paper, we write $k_i$ for the position of the $i$-th step, and \code{k0}, \code{k1}, \code{k2}, for the computer variables corresponding to ad hoc step positions depending on context.

The implementation closely mirrors the structure outlined in the pseudocode. It consists of a main function, \code{EscalierFit}, and two recursive helper functions, \code{buildESS} and \code{getTail}. The code requires $f(x)$ and its mean value $\bar f(a,b) = \frac{\int_{a}^{b} f(x) dx}{b-a}$ to be implemented. Additionally, the main function requires a precision argument for step position increments, a tolerance value for zero-testing, and an array pointer to return the  positions of optimal steps. For ease of implementation, the array is stored in reverse order. Unlike the pseudocode, the implemented algorithm does not directly return this array; instead, it modifies it internally and subsequently returns the ESS value.

To determine the optimal $k_1$ value, \code{EscalierFit} performs a sweep search with increments equal to the  precision argument. At each iteration of the main loop, it calls \code{buildESS} to compute the ESS of the current $k_1$ candidate, and keeps track of the running maximum (stored in \code{argmax} and \code{maxess} variables). After the loop is completed, \code{EscalierFit} supplies the optimal $k_1$ value to \code{getTail} to retrieve the array of optimal $k$ values.

Given arguments \code{k0} and \code{k1}, the helper functions \code{buildESS} and \code{getTail} aim to recursively maximize the ESS that can be achieved with several steps $k_0 < k_1 < k_2 < \dots < b$. The key distinction is that \code{buildESS} returns only the maximum ESS value, whereas \code{getTail} also stores the corresponding optimal $k$ values in the provided array pointer. Note that the sweep search for the optimal $k_1$ value in \code{EscalierFit} could be performed with the \code{getTail} function only, however, it is more efficient to use the \code{buildESS} function in terms of computational cost and memory management, then run the \code{getTail} function  once to recover the list of optimal $k$ values.

The pseudocode for both \code{buildESS} and \code{getTail} requires to solve an equation for the next critical $k$ value:
\[
    \text{Solve for }  k_{i+1} \in (k_i,b) :
    \displaystyle
    \frac{\int_{k_i}^{k_{i+1}} f(x) \d x}{k_{i+1}-k_i} = 2 f(k_i) - \frac{\int_{k_{i-1}}^{k_i} f(x) \d x}{k_i-k_{i-1}}
\]
This is implemented by another sweep search of all \code{k2} values that satisfy the equation up to the specified tolerance level.  Specifically, the `criticality' of a candidate \code{k2} value is calculated as
\[
 \norm*{\bar f(\code{k0}, \code{k1}) + \bar f(\code{k1},\code{k2}) - 2f(\code{k1})} = \norm*{\frac{\int_{\smallcode{k1}}^{\smallcode{k2}} f(x) \d x}{\code{k2}-\code{k1}} +  \frac{\int_{\smallcode{k0}}^{\smallcode{k1}} f(x) \d x}{\code{k1}-\code{k0}} - 2f(\code{k1})}
\]
A criticality value of 0 represents perfect equality up to machine precision, so values closer to 0 indicate configurations more likely to produce optimal results. Listing~\ref{lst:czone} shows the corresponding code implementation together with the critical zone logic described below.

As noted earlier, the pseudocode for \code{buildESS} requires all $k_{i+1}$ solutions to be solved for, after which a separate recursive call is invoked for each. In order to implement that, \code{buildESS} loops through all possible \code{k2} values between \code{k1} and \code{b} in increments equal to the specified precision, and computes the criticality \code{c} of each \code{k2}. If \code{c} is less than the tolerance threshold, the code enters a \emph{critical zone} which corresponds to a secondary loop over subsequent \code{k2} values whose criticality \code{czone} remain below the threshold.  Inside a critical zone, \code{buildESS} tracks the \code{k2} value that has the lowest criticality using the \code{argmin} and \code{c} variables. As soon as the secondary loop exits the critical zone, \code{buildESS} recurses with one fewer steps with new arguments $\code{k0} \leftarrow \code{k1}$ and $\code{k1} \leftarrow \code{argmin}$. 

\lstdefinestyle{snippetStyle}{
    backgroundcolor=\color{backcolor},   
    commentstyle=\color{codegreen},
    keywordstyle=\color{blue},     numberstyle=\tiny\color{codegrey},
    stringstyle=\color{codepurple},
    basicstyle=\ttfamily\footnotesize,
    escapeinside={(*@}{@*)},
    breakatwhitespace=false,         
    breaklines=false,                 
    captionpos=b,                    
    keepspaces=true,                 
    numbers=none,                                     
    showspaces=false,                
    showstringspaces=false,
    showtabs=false,                  
    tabsize=3,
}

\begin{lstlisting}[
                   label={lst:czone},
                   caption ={Critical Zone Code Implementation},
                   style=snippetStyle,
                   float=h]
    c = abs(fbar(k0, k1) + fbar(k1, k2) - 2 * f(k1));
    if (c > EPSILON && c < tolerance) { // critical zone
		argmin = k2;
		do { 
			i++; 
			if ((k2 = k1 + i * precision) >= b)
				break;
			czone = abs(fbar(k0, k1) + fbar(k1, k2) - 2 * f(k1));
			if (c > czone) {
				c = czone;
				argmin = k2;
			}
		} while (czone < tolerance);
		maxess = max(maxess, buildESS(k1, argmin, b, maxsteps - 1));
	}
\end{lstlisting}

In practice, the critical zone approach produces excellent results for many of our test functions. The main limitation is that this logic assumes that each critical zone contains only one critical \code{k2} value, which is assumed to have the smallest criticality. This assumption may not always hold mathematically. Here, the tolerance parameter plays a crucial role in shaping the results: a lower tolerance yields smaller critical zones and may overlook potential solutions, while a higher tolerance broadens the zone but still restricts selection to a single candidate, potentially discarding other critical \code{k2} values.  Hence, there is a ``sweet spot'' for the tolerance parameter which ultimately depends on the curve $f(x)$ being fitted, the precision parameter and the maximum number of steps.

\subsection{Group I results} In this subsection, we examine Escalier’s behavior when applied to discontinuous functions (``Group I functions''): $f_2$, a non-square-integrable function with midpoint singularity; $f_4$, which is a unit step function with a step at $x=1$; $f_5$, a square-integrable function with midpoint singularity; and $f_8$,  a piecewise linear function with a jump discontinuity at $x=1$. Since the critical point equations \eqref{eq:critical-points-1} require $f$ to be continuous, Escalier is expected to fail for all Group I functions. For this group, our testing was performed using tolerance levels \(10^{-4}, 10^{-3}, 10^{-2}, 10^{-1}, 1\), precision levels \(10^{-4}, 10^{-3}, 10^{-2}\), and maximum step counts ($N_{\max}$) of 1 to 10.

To begin, we examined the performance of Escalier on the simplest case, $f_4$, which may serve as a sanity check. Since $f_4$ is a unit step function and thus discontinuous, we may wonder whether Escalier is able to approximate it given reasonable parameters. As it turns out, across all tolerance levels, precision settings, and step counts tested, Escalier consistently returned a perfect match. Notably, even when given $N_{\max} > 1$, the algorithm correctly identified a single step at $x=1$. To confirm that this behavior was not specific to $f_4$, we tested Escalier on thousands of randomly generated single-step functions, with step positions uniformly sampled from \((a,b)\) and heights from \((-10,10)\); in every case, Escalier produced a perfect match.

Reflecting on our implementation, Escalier employs a brute-force approach to determine the initial step, $k_1$, after which it utilizes the critical zone equation to narrow down the search for subsequent steps. Hence, there is no need for continuity when working with a single step. This clarifies why Escalier consistently achieved perfect approximations for single-step functions, such as $f_4$ and the randomly generated single-step tests. This observation was further confirmed by testing Escalier with all functions listed in Table~\ref{tab:testfunctions}, using $N_{\max}=1$, and we found that the algorithm consistently matched Brute Force results for all precision and tolerance settings.

\begin{figure}[htbp]
    \centering
    \begin{subfigure}[b]{0.49\textwidth}
        \includegraphics[width=\textwidth]{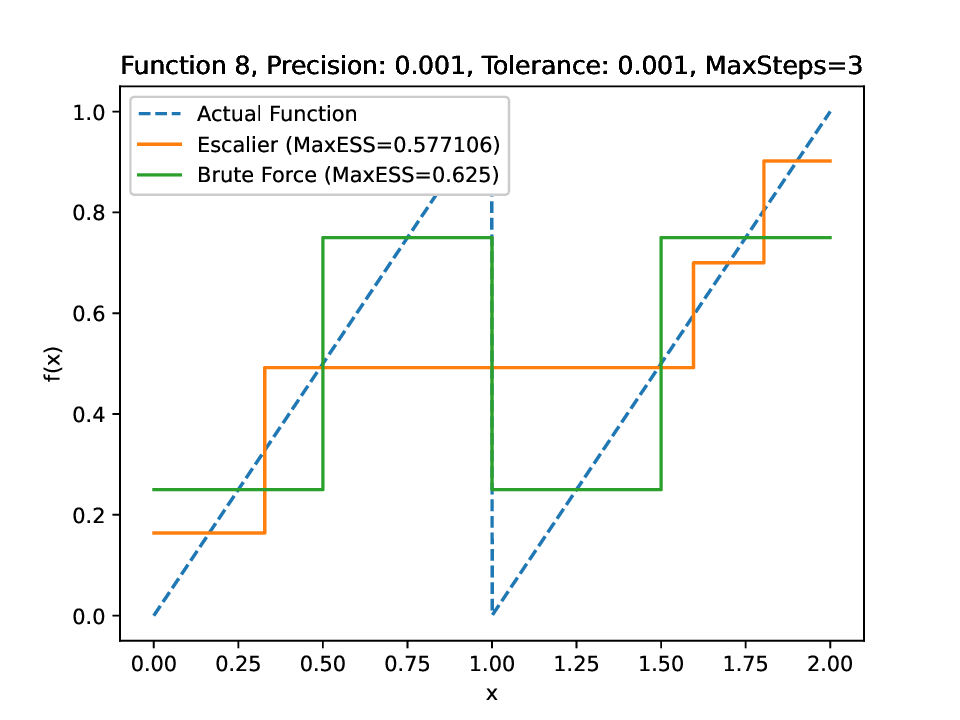}
        \caption{$N_{\max}=3$}
        \label{fig:f8subA}
    \end{subfigure}
    \hfill
    \begin{subfigure}[b]{0.49\textwidth}
        \includegraphics[width=\textwidth]{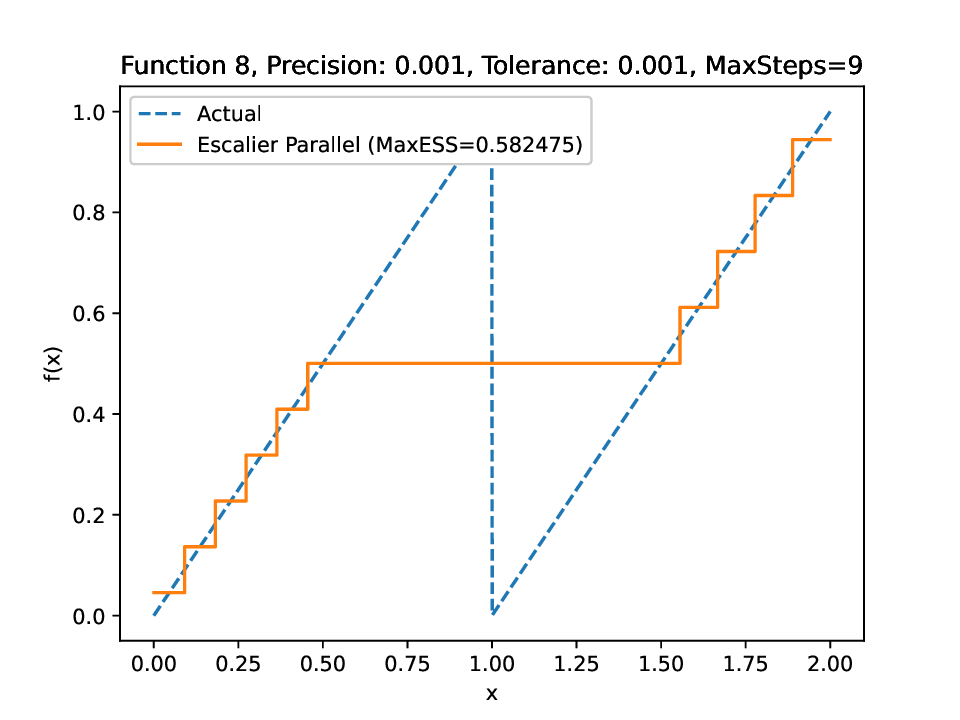}
        \caption{$N_{\max}=9$}
        \label{fig:f8subB}
    \end{subfigure}
    \hfill
    
    \caption{Approximation plot of $f_8$ with tolerance and precision = $10^{-3}$}
    \label{fig:f8PrecisionGraphs}
\end{figure}

Limitations begin to appear when Escalier is tested against a broader range of discontinuous functions. For randomly generated two-step functions, Escalier consistently failed to yield accurate results.  This is attributable to the fact that the algorithm relies on the critical zone equation to determine the optimal $k_2,k_3, \dots$ positions, which requires continuity.
We further investigated Escalier's performance on $f_8$, a double ramp function with discontinuity at $x=1$. It is easy to see that, for an odd number of steps, the optimal solution is made of regular steps across both ramps, with the middle step starting precisely at the jump. This optimal solution was correctly found by Brute Force tests for 3 steps (precision $10^{-3}$, see Figure~\ref{fig:f8subA}) and 5 steps (precision $10^{-2}$). In all our tests, Escalier failed to produce correct results, as illustrated in Figure~\ref{fig:f8subB} for tolerance and precision of $10^{-3}$, and $N_{\max}=9$.

Both $f_2$ and $f_5$ have midpoint singularities at $x=1$, but they differ in integrability: $f_2$ is non-square-integrable, therefore it lies outside the $L^2$ Hilbert space of our theoretical framework (Sections \ref{sec:intro} to \ref{sec:nonlinear-lsq}), whereas $f_5$ is square-integrable. Despite this distinction, Escalier behaves similarly for both functions: for all parameter choices, it placed a single step at the singularity and returned a constant ESS value of 8 for $f_2$, and 18 for $f_5$.

In summary, in agreement with theory, Escalier usually fails to deliver meaningful results in the presence of discontinuities and singularities. However, Escalier produces perfect approximations for a single unit step function, and more generally identifies a single discontinuity when using $N_{\max}=1$.

\subsection{Group II results} In this subsection, we evaluate the performance of Escalier on Group II test functions, which are both square-integrable and continuous on $(a,b)$: $f_1$, a smooth and monotonic function; $f_3$, a square-integrable function with an endpoint singularity; $f_6$, a continuous, one-dimensional L\'evy function; and $f_7$, a smooth oscillatory function.

As previously mentioned, performance was evaluated in terms of speed and $R^2$. To reduce any noise in execution time, all runtimes were measured as the average of three successive runs. Additionally, parallel versions of both Escalier and Brute Force were implemented and tested alongside their serial counterparts, providing a broader range of results for analysis. Testing was carried out as follows:
\begin{enumerate}[leftmargin=*]
    \item Baseline comparison:  We compared wall clock runtime and $R^2$ accuracy of Escalier and Brute Force algorithms for $N_{\max}=2$ and $N_{\max}=3$. 
    \item Scalability analysis: We examined the increase in Escalier's runtime and accuracy as $N_{\max}$ varies from 2 to 10.
\end{enumerate}

\subsubsection{Baseline comparison: Escalier versus Brute-Force ($N_{\max}=2,3$)} We tested Group II test functions using tolerance levels $10^{-4}, 10^{-3}, 10^{-2}, 10^{-1}, 1$, and  precision levels $10^{-4}, 10^{-3},$ and $10^{-2}$. We begin by comparing the runtime efficiency of both algorithms, where Table~\ref{tab:Runtime-5.3.1} reports the runtime (in milliseconds) for their serial and parallel implementations. Note that Table~\ref{tab:Runtime-5.3.1-B} excludes the serial runtimes for precision $10^{-4}$, as Brute Force computation time was otherwise prohibitive.

\begin{figure}[htbp]            
    \centering
    \includegraphics[width=\linewidth]{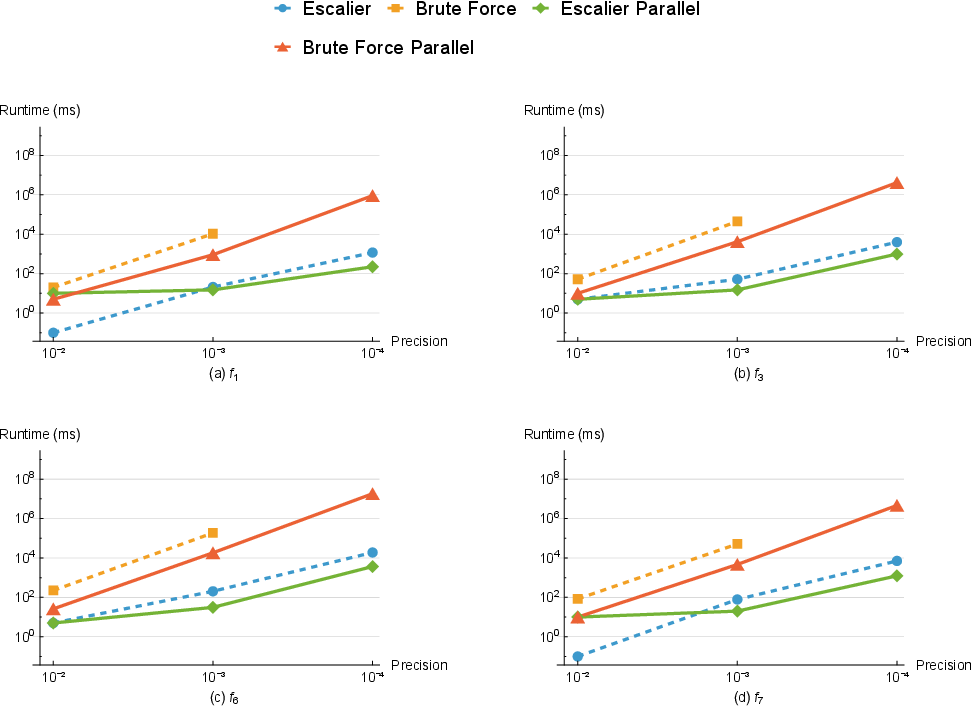} 
    \vspace{0.5em}
    \caption{Log–log plots of runtime (ms) against precision for Escalier and Brute Force algorithms, $N_{\max}=3$, tolerance =$10^{-3}$}
    \label{fig:runtime-vs-precision}
\end{figure}

For $N_{\max}=2$ and precision $10^{-2}$, runtimes across all functions are nearly identical. This extends to the parallel implementations at all precision levels for this step count, where differences remain negligible. In contrast, the serial implementations begin to diverge as the precision becomes finer, though the gap is not yet drastic. This is expected given Escalier's implementation which involves an exhaustive sweep search over all $k_1$ values, while subsequent $k_i$ searches are narrowed down within the critical zone. The impact becomes more pronounced for $N_{\max}=3$. With the exception of parallel implementations at precision $10^{-2}$, Escalier shows a clear runtime advantage over Brute-Force across all precision levels, with the differences widening substantially at finer precisions. Figure~\ref{fig:runtime-vs-precision}, which corresponds to runtimes from Table~\ref{tab:Runtime-5.3.1-B} along with additional entries for serial Escalier at precision $10^{-4}$, illustrates these trends. The log–log plots show the runtime growth for all functions and algorithms as the precision becomes finer, confirming that Escalier consistently scales more favorably than Brute Force, with negligible differences only at precision $10^{-2}$. Additionally, we observe that the complexity of the function influences runtime performance. Escalier achieves its lowest runtimes with simpler, monotone functions such as $f_1$, whereas more complex, highly oscillatory functions such as $f_6$ require greater computational efforts.

\begin{table}[htbp]
    \centering
    \caption{Runtime (ms) for Escalier (ES) and Brute‑Force (BF) algorithms for $N_{\max} = 2,3$ at various precision parameters.  Escalier tolerance = $10^{-3}$. Parallel implementation is indicated with $\varparallel$ symbol.}
    \label{tab:Runtime-5.3.1}

    \begin{subtable}{\linewidth}
        \centering
        \caption{\(N_{\max}=2\)}
        \begin{tabular}{c   *{4}{T{2}} *{4}{T{3}}  *{2}{T{5}} *{2}{T{4}}   }
        \toprule
                        & \multicolumn{4}{c}{Precision $=10^{-2}$}
            & \multicolumn{4}{c}{Precision $=10^{-3}$}
            & \multicolumn{4}{c}{Precision $=10^{-4}$}\\
        \cmidrule(lr){2-5} \cmidrule(lr){6-9} \cmidrule(lr){10-13}
                & {ES} & {BF} & {ES}$_{\varparallel}$ & {BF}$_{\varparallel}$ 
        & {ES} & {BF} & {ES}$_{\varparallel}$ & {BF}$_{\varparallel}$ 
        & {ES} & {BF} & {ES}$_{\varparallel}$ & {BF}$_{\varparallel}$  \\
        \midrule
        $f_1$ & 0  & 0 & 10 & 5 & 15  & 41  & 10 & 5  & 750   & 1328  & 93   & 120\\
        $f_3$ & 5  & 5 & 5  & 0 & 67  & 93  & 21 & 10 & 3307  & 4838  & 448  & 458\\
        $f_6$ & 10 & 5 & 10 & 0 & 151 & 250 & 36 & 31 & 14146 & 24599 & 1865 & 1994\\
        $f_7$ & 0  & 5 & 5  & 5 & 57  & 94  & 15 & 15 & 3921  & 6880  & 479  & 531\\
        \bottomrule
        \end{tabular}
    \end{subtable}
    
    \vspace{1.0ex}

    \begin{subtable}{\linewidth}
        \centering
        \caption{\(N_{\max}=3\)}
        \begin{tabular}{c   T{2}T{3}T{2}T{2}   T{3}T{6}T{2}T{5}   T{4}T{8}}
            \toprule
                                & \multicolumn{4}{c}{Precision $=10^{-2}$}
                & \multicolumn{4}{c}{Precision $=10^{-3}$}
                & \multicolumn{2}{c}{Precision $=10^{-4}$}\\
            \cmidrule(lr){2-5} \cmidrule(lr){6-9} \cmidrule(lr){10-11}
                        & {ES} & {BF} & {ES}$_{\varparallel}$ & {BF}$_{\varparallel}$ 
            & {ES} & {BF} & {ES}$_{\varparallel}$ & {BF}$_{\varparallel}$ 
            & {ES}$_{\varparallel}$ & {BF}$_{\varparallel}$  \\
            \midrule
            $f_1$ & 0 & 20  & 10 & 5  & 21  & 10573  & 15 & 911   & 224  & 907760 \\
            $f_3$ & 5 & 52  & 5  & 10 & 52  & 44734  & 15 & 4161  & 989  & 4150276 \\
            $f_6$ & 5 & 224 & 5  & 26 & 203 & 186458 & 31 & 18245 & 3666 & 18167359 \\
            $f_7$ & 0 & 83  & 10 & 10 & 78  & 51343  & 20 & 4661  & 1239 & 4624547 \\
        \bottomrule
        \label{tab:Runtime-5.3.1-B}
    \end{tabular}
    \end{subtable}

    \vspace{0.5ex}
        \end{table}

Table~\ref{tab:ESS-comparison-N2} reports the R\textsuperscript{2} values obtained by Escalier and Brute Force for Group II functions at each precision level with $N_{\max}=2$, while Table~\ref{tab:ESS-comparison-N3} reports the corresponding results for $N_{\max}=3$.  It is worth noting that, while Brute Force serves as our baseline for the best achievable escalier fit, it may not match the theoretical optimum given by \eqref{eq:critical-points-2} that Escalier's logic is based upon.  As such, comparing the $R^2$ of Escalier against Brute Force may not be an entirely accurate performance metric, particularly at cruder precision levels.

For $f_1$, a smooth monotonic function, Escalier produces exact matches at both \(N_{\max}=2\) and \(N_{\max}=3\) for precisions \(10^{-4}\) and \(10^{-3}\). Minor discrepancies appear only when the precision is relaxed to \(10^{-2}\). This suggests that the algorithm is robust when applied to simple, well-behaved functions.  \( f_3 \) is also smooth and monotonic, but introduces an endpoint singularity at \( x=0 \). Despite this potential challenge, Escalier performs similarly, producing exact matches at both \(N_{\max}=2\) and \(N_{\max}=3\) for precisions \(10^{-4}\) and \(10^{-3}\). This suggests that endpoint singularities may have a limited influence on the algorithm’s performance.

In contrast, \( f_6 \), which exhibits rapid oscillations, proves more challenging for Escalier. Although the algorithm achieves a reasonable $R^2$ with \(N_{\max}=2\) at finer precision levels ($10^{-4}$ and $10^{-3}$), more noticeable discrepancies appear at \(N_{\max}=3\), especially at $10^{-2}$ precision. While no parameter configuration produced a perfect match against Brute Force, Escalier often performed well in terms of $R^2$.  Further investigation with a broader range of parameters and tolerance settings may provide a deeper understanding of Escalier's limitations in such settings.

Finally, \( f_7 \) is nonmonotonic but exhibits fewer oscillations. For $N_{\max}= 2$, Escalier produces close approximations at precision levels $10^{-2}$ and $10^{-3}$, and achieves a perfect match with Brute Force at precision $10^{-4}$. For $N_{\max}= 3$, there is a slightly larger gap in accuracy at precision $10^{-2}$, yet the algorithm still achieves a near-perfect and a perfect approximation for precisions $10^{-3}$ and $10^{-4}$, respectively. This behavior suggests that Escalier performs well on functions that exhibit mild oscillations.

\begin{figure}[htbp]
  \centering
  \begin{subfigure}[b]{0.49\textwidth}
    \includegraphics[width=\textwidth]{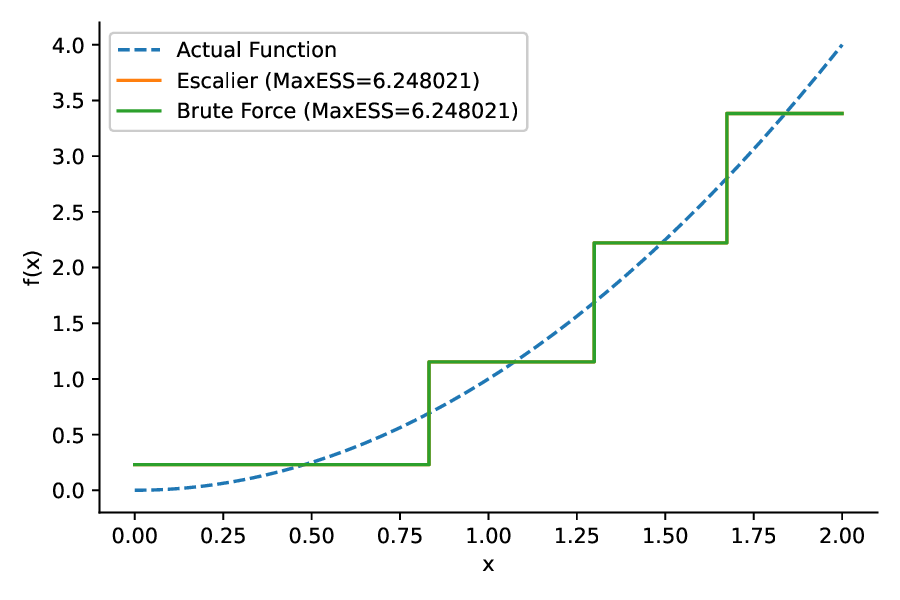}
    \caption{\(f_1(x)\)}
    \label{fig:subA}
  \end{subfigure}
  \hfill
  \begin{subfigure}[b]{0.49\textwidth}
    \includegraphics[width=\textwidth]{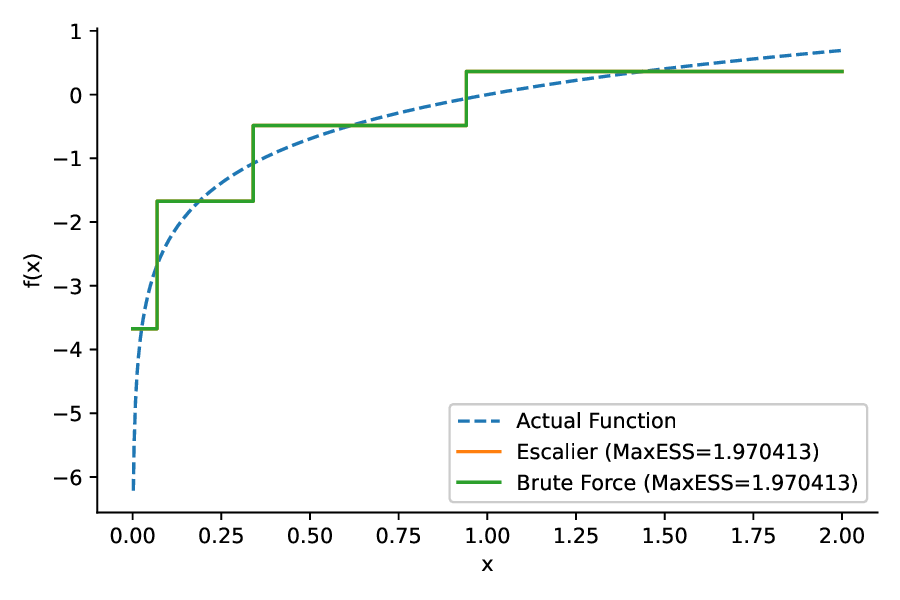}
    \caption{\(f_3(x)\) }
    \label{fig:subB}
  \end{subfigure}
  
  \vskip\baselineskip
  
  \begin{subfigure}[b]{0.49\textwidth}
    \includegraphics[width=\textwidth]{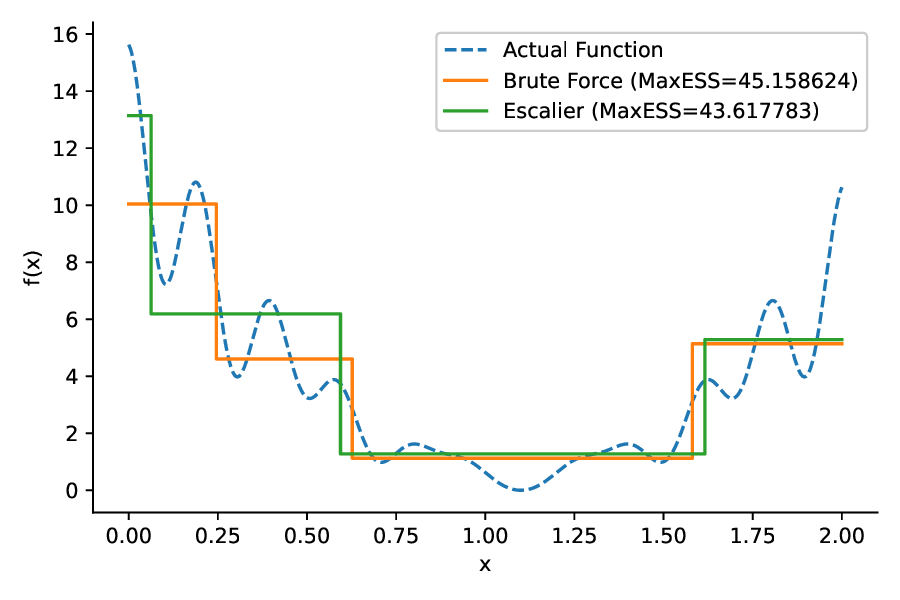}
    \caption{\(f_6(x)\)}
    \label{fig:subC}
  \end{subfigure}
  \hfill
  \begin{subfigure}[b]{0.49\textwidth}
    \includegraphics[width=\textwidth]{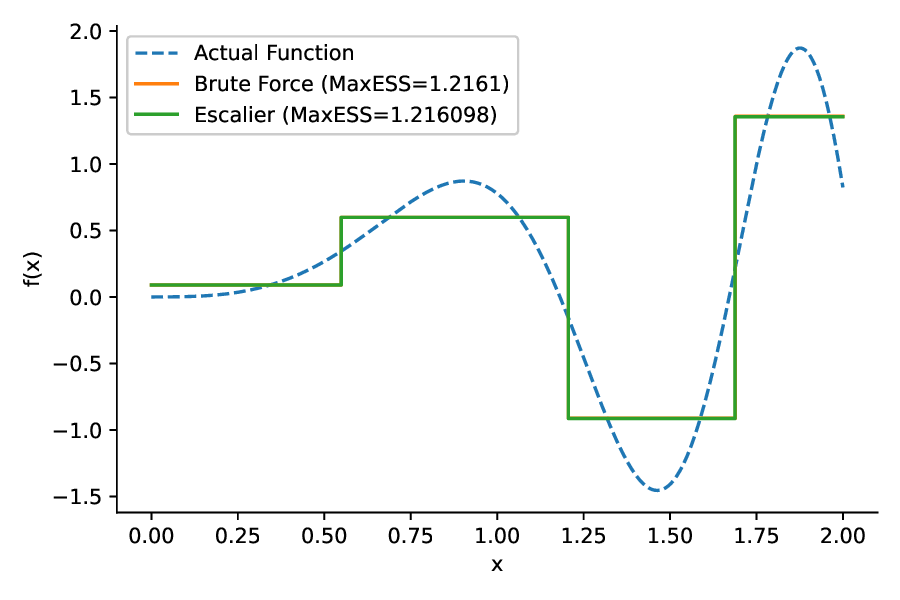}
    \caption{\(f_7(x)\)}
    \label{fig:subD}
  \end{subfigure}
  
  \caption{Group II graph comparison for default tolerance and precision of $10^{-3}$ and $N_{\max}=3$.}
  \label{fig:comparison}
\end{figure}

Additionally, we observe a clear relationship between precision and Escalier's approximation quality as measured by $R^2$. Across all functions and both step counts $N_{\max}=2,3$, finer precision improved Escalier's $R^2$. When these results are compared with our runtime data in Table~\ref{tab:Runtime-5.3.1}, the additional computation time required by finer precision appears worthwhile.  While for Brute Force finer precisions may only yield a modest gain in accuracy, for Escalier this can mean a jump in approximation quality at a mere fraction of the Brute Force runtime. The default case at $N_{\max}=3$  (with tolerance and precision of $10^{-3}$) is shown in Figure ~\ref{fig:comparison} as a representative example.

\begin{table}[htbp]
    \small
    \centering
    \caption{R\textsuperscript{2} comparison of Escalier vs.\ Brute-Force }
    \label{tab:ESS-comparison-5.3.1}

    \begin{subtable}{\linewidth}
    \centering
    \caption{\(N_{\max}=2\)}
    \label{tab:ESS-comparison-N2}
    \begin{tabular}{c ccc ccc ccc}
        \toprule
                        & \multicolumn{3}{c}{Precision $=10^{-2}$}
            & \multicolumn{3}{c}{Precision $=10^{-3}$}
            & \multicolumn{3}{c}{Precision $=10^{-4}$}\\
        \cmidrule(lr){2-4} \cmidrule(lr){5-7} \cmidrule(lr){8-10}
                & {Escalier} & {BF} & $\Delta$
        & {Escalier} & {BF} & $\Delta$
        & {Escalier} & {BF} & $\Delta$\\
        \midrule
        $f_1$ & 0.956888 & 0.957022 & 0.000134 & 0.957025 & 0.957025 & 0 & 0.957025 & 0.957025 & 0 \\
        $f_3$ & 0.833849 & 0.835712 & 0.001863 & 0.835731 & 0.835731 & 0 & 0.835731 & 0.835731 & 0 \\
        $f_6$ & 0.799839 & 0.867428 & 0.067589 & 0.866807 & 0.867498 & 0.000691 & 0.867440 & 0.867499 & 0.000059 \\
        $f_7$ & 0.793288 & 0.797437 & 0.004149 & 0.797610 & 0.797669 & 0.000059 & 0.797670 & 0.797670 & 0 \\
        \bottomrule
    \end{tabular}
    \end{subtable}
    
    \vspace{1.0ex}

    \begin{subtable}{\linewidth}
    \centering
    \caption{\(N_{\max}=3\)}
    \label{tab:ESS-comparison-N3}
    \begin{tabular}{c ccc ccc ccc}
        \toprule
                        & \multicolumn{3}{c}{Precision $=10^{-2}$}
            & \multicolumn{3}{c}{Precision $=10^{-3}$}
            & \multicolumn{3}{c}{Precision $=10^{-4}$}\\
        \cmidrule(lr){2-4} \cmidrule(lr){5-7} \cmidrule(lr){8-10}
                & {Escalier} & {BF} & $\Delta$
        & {Escalier} & {BF} & $\Delta$
        & {Escalier} & {BF} & $\Delta$\\
        \midrule
        $f_1$ & 0.976133 & 0.976246 & 0.000113 & 0.976253 & 0.976253 & 0 & 0.976253 & 0.976253 & 0\\
        $f_3$ & 0.865532 & 0.900420 & 0.034888 & 0.900424 & 0.900424 & 0 & 0.900424 & 0.900424 & 0\\
        $f_6$ & 0.799839 & 0.926175 & 0.126336 & 0.894766 & 0.926374 & 0.031608 & 0.924830 & 0.926376 & 0.001546\\
        $f_7$ & 0.793288 & 0.849557 & 0.056269 & 0.849694 & 0.849696 & 0.000002 & 0.849699 & 0.849699 & 0\\
        \bottomrule
    \end{tabular}
    \end{subtable}

    \vspace{0.5ex}
\end{table}

All previous Escalier results were obtained with standard tolerance parameter of $10^{-3}$. We now examine how changing tolerance affects Escalier's accuracy. As mentioned earlier, tolerance directly governs the solution selection process: lower tolerance yields smaller critical zones, potentially skipping viable solutions, while higher tolerance broadens the critical zone but potentially disregards some valid \code{k2} values within. Through systematic experimentation, we find that there is no ``one-size-fits-all'' optimal tolerance value. Rather, the appropriate tolerance parameter depends on the precision level, maximum step count, and function $f$. We found that tolerance can help improving Escalier's performance with more complex functions like $f_6$ and $f_7$.

\begin{figure}[htbp]
    \centering
    \begin{subfigure}[b]{0.49\textwidth}
    \includegraphics[width=\textwidth]{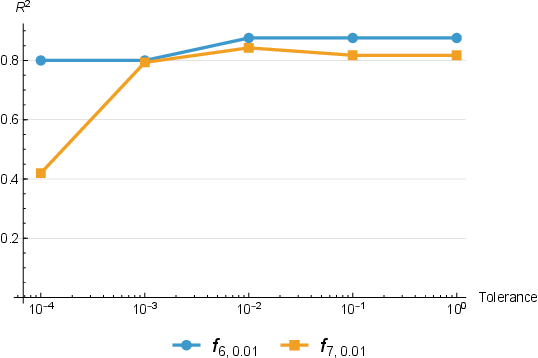}
    \caption{Precision = \(10^{-2}\)}
    \label{fig:ToleranceVsR2-3-p0.01}
  \end{subfigure}
  \hfill
  \begin{subfigure}[b]{0.49\textwidth}
    \includegraphics[width=\textwidth]{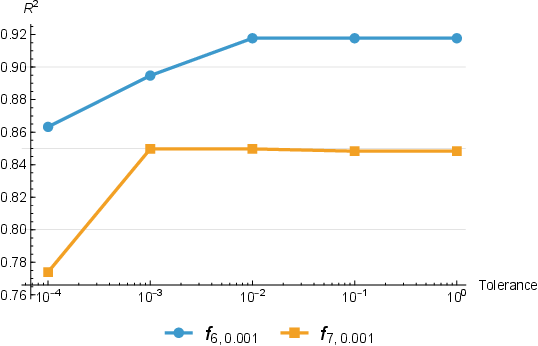}
    \caption{Precision = \(10^{-3}\)}
    \label{fig:ToleranceVsR2-3-p0.001}
  \end{subfigure}
  \hfill
    
\vskip\baselineskip

  \begin{subfigure}[b]{0.49\textwidth}
    \includegraphics[width=\textwidth]{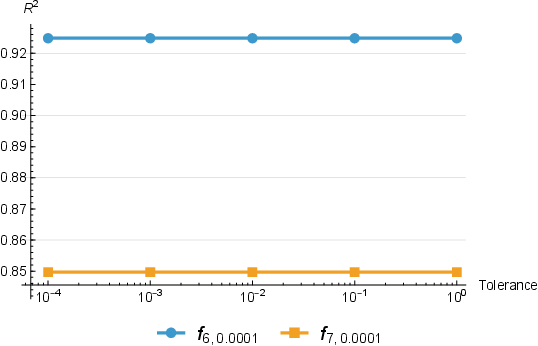}
    \caption{Precision = \(10^{-4}\)}
    \label{fig:ToleranceVsR2-3-p0.0001}
  \end{subfigure}
  
    \caption{Tolerance vs R\textsuperscript{2} for $f_6$ and $f_7$ with $N_{\max} =3$ }
    \label{fig:ToleranceVsR2-3}
\end{figure}

\begin{figure}[htbp]
  \centering
    \begin{subfigure}[b]{0.49\linewidth}
    \includegraphics[width=\linewidth]{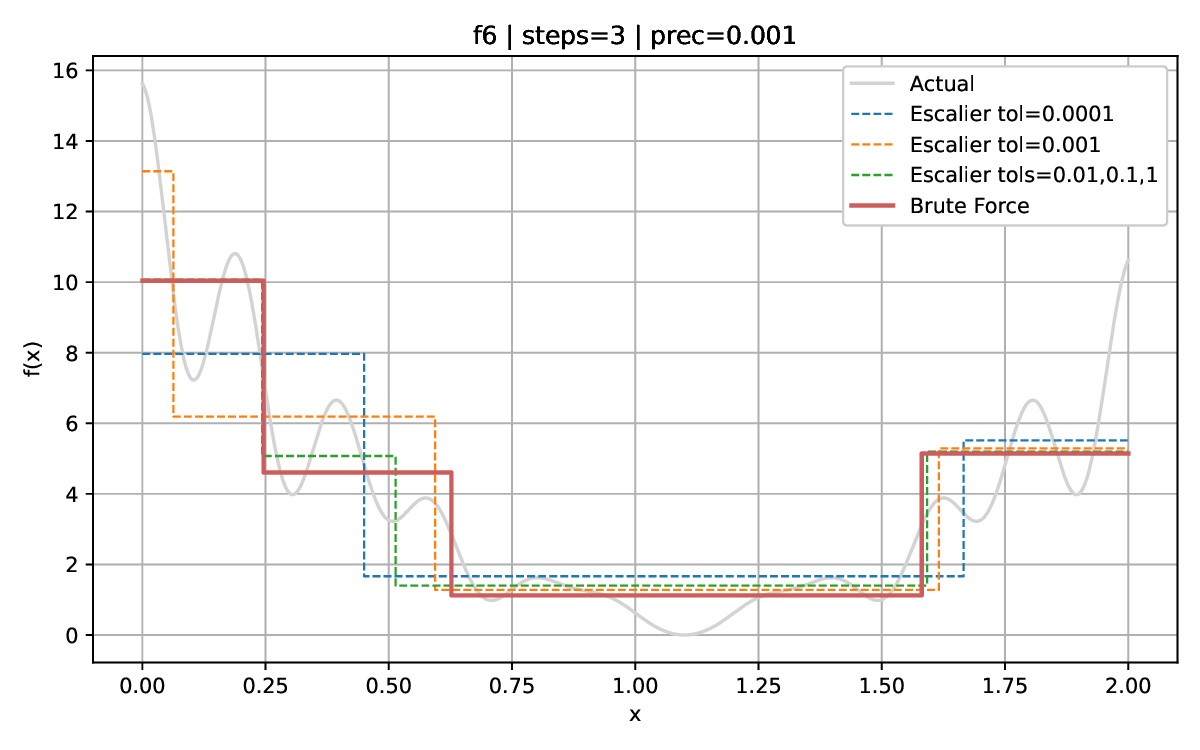}
    \caption{$f_6$}
    \label{fig:subA}
  \end{subfigure}\hfill
    \begin{subfigure}[b]{0.49\linewidth}
    \includegraphics[width=\linewidth]{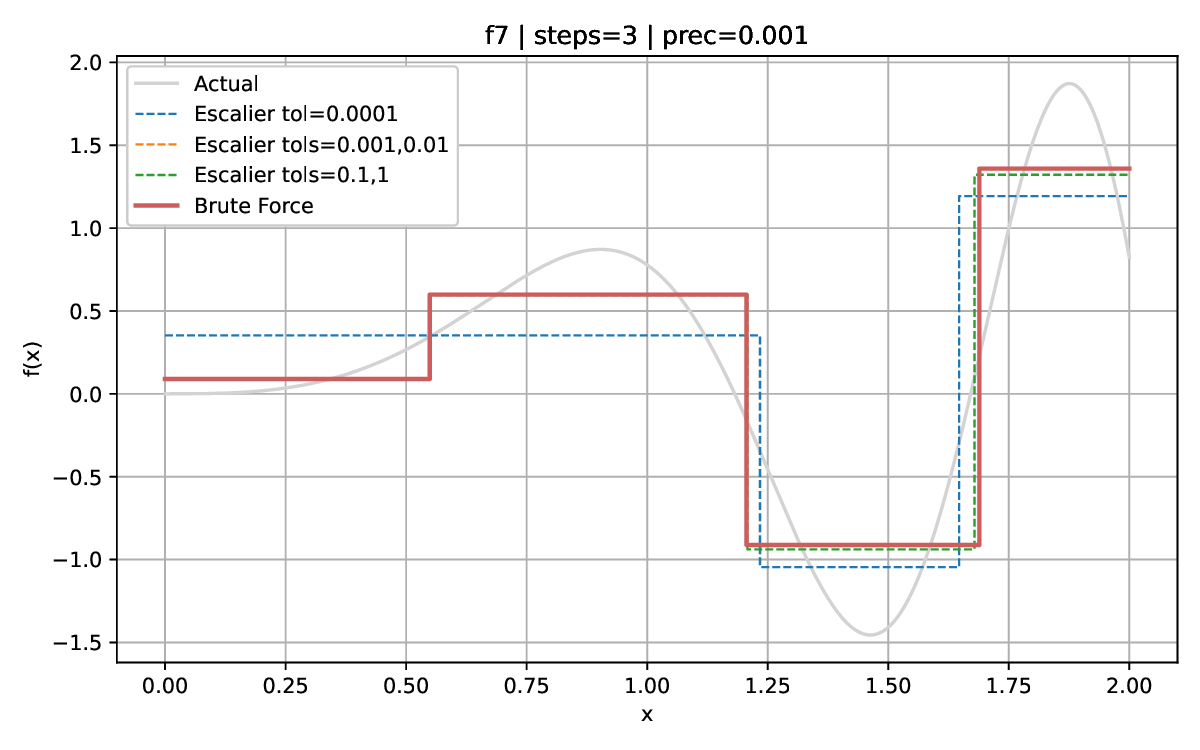}
    \caption{$f_7$}
    \label{fig:subB}
  \end{subfigure}

  \caption{Impact of Tolerance Variation on Approximation quality, for Precision = $10^{-3}$; $N_{\max}=3$.}
  \label{fig:ToleranceVariationPlots}
\end{figure}

Figures~\ref{fig:ToleranceVsR2-3} and \ref{fig:ToleranceVariationPlots} highlight how changing the tolerance level, namely, $10^{-4}$, $10^{-3}$, $10^{-2}$, $10^{-1}$, $1$ affects the $R^2$ accuracy of the approximation for $f_6$ and $f_7$, for various precision levels. Results for $f_1$ and $f_3$ are not shown, as the monotone functions already achieve high accuracy, making tolerance effects minimal. For $f_6$, the $R^2$-tolerance curve is nondecreasing for all precisions: at $10^{-2}$ and $10^{-3}$, the $R^2$ rises then plateaus at tolerance $10^{-2}$, with no further gain beyond that point (see Figures~\ref{fig:ToleranceVsR2-3-p0.01} and ~\ref{fig:ToleranceVsR2-3-p0.001}); at precision $10^{-4}$, Escalier yields the same $R^2$ across all tolerance levels tested (see Figure~\ref{fig:ToleranceVsR2-3-p0.0001}). This was not the case for $f_7$, where at precision $10^{-2}$,  a clear peak in $R^2$ occurs at a tolerance of $10^{-2}$. As the precision gets finer, this peak becomes less pronounced, with identical optimal $R^2$ achieved at tolerance levels $10^{-3}$ and $10^{-2}$. At precision $10^{-4}$, the behavior resembles that of $f_6$, where tolerance variation does not have an effect on $R^2$. These findings underscore  that the optimal tolerance level, ``sweet spot'', varies according to the structure of the function and the chosen precision level.

\subsubsection{Scalability analysis $(N_{\max}\geq 2$)} In this subsection, we investigate the scalability of Escalier by examining the increase in runtime and $R^2$ as the maximum step count $N_{\max}$ varies from 2 to 10. During this analysis, we encountered two main difficulties. First, under the default settings (tolerance and precision of $10^{-3}$), the parallel implementation was so efficient that all runtimes were nearly instant, making it difficult to observe growth. Second, identifying a single tolerance to represent growth proved challenging. In some scenarios, certain tolerance levels led the algorithm to plateau, making it difficult to capture the true scaling behavior. This issue is discussed in greater detail later in the subsection. To resolve these difficulties, we adopted a finer precision of $10^{-5}$. This adjustment amplifies runtime growth across step counts, thereby making scaling behavior observable. Moreover, as discussed previously, employing a finer precision also reduces variability in results due to tolerance, providing a more stable setting in which to evaluate scalability.

\begin{figure}[htbp]
  \centering

  \begin{subfigure}[b]{0.49\linewidth}
    \includegraphics[width=\linewidth]{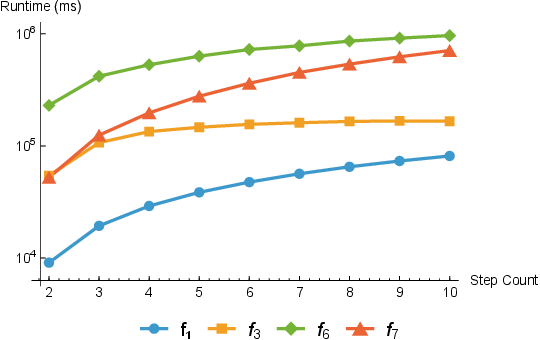}
    \caption{Runtime}
    \label{fig:Runtime-5.3.2}
  \end{subfigure}\hfill
    \begin{subfigure}[b]{0.49\linewidth}
    \includegraphics[width=\linewidth]{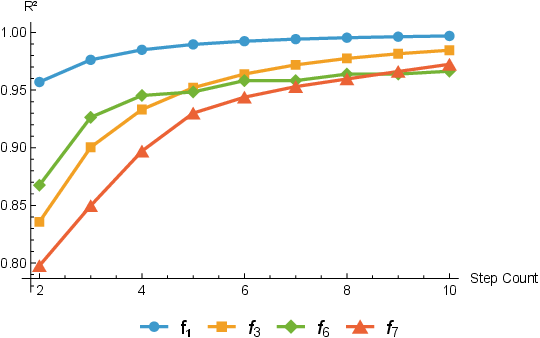}
    \caption{$R^2$}
    \label{fig:R2-5.3.2}
  \end{subfigure}

  \caption{Escalier Parallel's runtime and $R^2$ for $N_{\max}= 2-10$ , tolerance $= 10^{-3}$, and precision $= 10^{-5}$}
  \label{fig:5.3.2-Steps-time-R2-plots}
\end{figure}

Table~\ref{tab:Runtime-R2-5.3.2} reports the runtimes and $R^2$ values  of Escalier Parallel as $N_{\max}$ increases from 2 to 10, with tolerance $10^{-3}$ and precision $10^{-5}$. As observed previously, the complexity of the test function strongly governs runtime, with $f_{6}$ consistently the most demanding, followed by $f_{7}$, $f_{3}$, and $f_{1}$. At this precision level, runtimes become significant, reaching approximately $9.6\times 10^{5}$ ms (about 16 minutes) for $f_{6}$ at $N_{\max}=10$, but they remain feasible for large-scale testing. Figure \ref{fig:Runtime-5.3.2} reproduces the runtime data on a logarithmic scale. The curves show clear growth as step count increases, yet the trajectories are sub-exponential. The slope flattens as $N_{\max}$ grows, most visibly for $f_{3}$, where the incremental cost per step decreases steadily. This indicates growth closer to linear than exponential across all tested functions.

Across all four test functions, $R^2$ improves steadily as more steps are introduced, with the exception of $f_{6}$ at $N_{\max}=7$ and $N_{\max}=9$, where the values remain unchanged. These plateaus highlight occasional irregularities in how tolerance interacts with more complex functions. Figure \ref{fig:R2-5.3.2} illustrates the overall trend. The steepest gains occur at small step counts, especially between $N_{\max}=2$ and $N_{\max}=5$, after which improvements taper off. This diminishing return indicates that the bulk of approximation accuracy is achieved early, while later steps contribute progressively smaller refinements. 

\begin{table}[htbp]
    \centering
    \caption{Escalier parallel runtime and $R^2$ for $N_{\max}= 2,3,\dots,10$, tolerance $= 10^{-3}$, and precision $= 10^{-5}$.}
    \label{tab:Runtime-R2-5.3.2}

        \begin{subtable}{\linewidth}
        \centering
        \caption{Runtime (ms)}
        \label{tab:Runtime-5.3.2}
        \begin{tabular}{c  *{9}{T{6}}}
            \toprule
                                                            $N_{\max}$ & {2}& {3}& {4}& {5}& {6}& {7}& {8}& {9}& {10} \\
            \midrule
            $f_1$ & 9104 & 19395 & 29161 & 38588 & 47583 & 56468 & 65177 & 73406 & 81380 \\
            $f_3$ & 54333 & 106948 & 134286 & 146682 & 155822 & 160948 & 165348 & 166932 & 165947 \\
            $f_6$ & 229932 & 418343 & 529135 & 631093 & 723677 & 780640 & 859968 & 914062 & 963135 \\
            $f_7$ & 52375 & 124718 & 197135 & 277369 & 361979 & 451489 & 535849 & 622677 & 708307 \\
        \bottomrule
    \end{tabular}
    \end{subtable}
  
    \vspace{1.0ex}

    \begin{subtable}{\linewidth}
        \centering
        \caption{R\textsuperscript{2}}
        \label{tab:R2-5.3.2}
        \begin{tabular}{c   *{9}{c} }
            \toprule
                                                           $N_{\max}$ & {2}& {3}& {4}& {5}& {6}& {7}& {8}& {9}& {10} \\
            \midrule
            $f_1$ & 0.957025 & 0.976253 & 0.984964 & 0.989633 & 0.992423 & 0.994221 & 0.995447 & 0.996321 & 0.996966 \\
            $f_3$ & 0.835731 & 0.900424 & 0.933200 & 0.952085 & 0.963955 & 0.971901 & 0.977480 & 0.981549 & 0.984606 \\
            $f_6$ & 0.867499 & 0.926376 & 0.945360 & 0.948399 & 0.958243 & 0.958243 & 0.963905 & 0.963905 & 0.966409 \\
            $f_7$ & 0.797670 & 0.849699 & 0.897083 & 0.930035 & 0.943853 & 0.953029 & 0.959662 & 0.966145 & 0.972407 \\
        \bottomrule
    \end{tabular}
    \end{subtable}

    \vspace{0.5ex}
\end{table}

As noted earlier, certain tolerance levels can cause the algorithm to plateau as the step count increases. To further investigate this point and clarify the relationship between tolerance and precision, Table~\ref{tab:R2-tolerance-vs-precision} presents the $R^2$ values for all Group II test functions. Each row corresponds to a function, while the two columns distinguish between precision levels $10^{-3}$ and $10^{-5}$. Within each graph, $R^2$ is shown as $N_{\max}$ increases from 2 to 10 at five tolerance levels, allowing direct comparison of how precision affects tolerance variation.

At precision $10^{-3}$ the effect of tolerance on $R^2$ is apparent. For the monotone functions $f_1$ and $f_3$, most tolerance levels yield consistent and steadily increasing $R^2$ values as the step count grows. For $f_1$, all tolerances produce identical results, with the exception of tolerance $10^{-4}$ where the curve plateaus at $N_{\max} = 3$. Test function $f_3$ follows a similar pattern, but a slight divergence among tolerance levels becomes visible when $N_{\max}\geq 8$. In contrast, $f_6$ exhibits far greater variability across its tolerance curves. Here, curves for tolerances $10^{-4}$ and $10^{-3}$ plateau almost immediately, at $N_{\max}=2$ and $N_{\max}=3$, respectively. The remaining tolerance curves continue to grow, but their relative performance shift with step count. For instance, at $N_{\max}=4$, a tolerance of $1$ ranks behind tolerances $10^{-1}$ and $10^{-2}$ , yet by $N_{\max}=10$ it surpasses both. This example supports our earlier claim that the ``best'' tolerance level is not fixed but depends jointly on the curve being approximated, the step count, and the precision setting. For $f_7$, the curve for tolerance $10^{-4}$ plateaus immediately at $N_{\max}=2$, while all other tolerance curves continue to improve as the step count increases. Among these, clear differences remain visible: as the step count grows, a tolerance of $10^{-2}$ consistently achieves the highest $R^2$, followed, in order, by $10^{-3}, 10^{-1}$, and $1$. 

Across all functions at this precision level, a tolerance of $10^{-4}$ consistently produces plateaus. Looking back at the critical zone implementation, it is clear that a small tolerance will narrow down our critical zone condition, which we can think of as our allowed error margin from the exact solution. From here, we may infer that precision $10^{-3}$ is too coarse to capture $k_{i+1}$ values that falls within such interval. By choosing a finer precision, we expand the search space and enable the algorithm to find more precise solutions that are more likely to satisfy the critical zone equation.

With finer precision at $10^{-5}$, the variation in $R^2$ due to tolerance is considerably reduced or eliminated for all functions. For $f_1$ and $f_3$, all tolerance levels produce identical results. Function $f_7$ follows the same trend, with only imperceptible differences: tolerance values of $10^{-1}$ and $1$ yield slightly lower $R^2$ values for $N_{\max}=8,9,10$, though this is not visible in the plots. Function $f_6$ shows minor but noticeable variation, yet far less than at coarse precision. Overall, $f_1$, $f_3$, and $f_7$ exhibit smooth growth as the step count increases, while $f_6$ also improves but in a less regular manner, which is expected given its complexity.

Overall, these observations confirm that the optimal tolerance level is not fixed, but depends on the function, the step count, and the precision setting. At the same time, the results demonstrate that $R^2$ variation due to tolerance can be effectively minimized, if not eliminated, by employing finer precision, at the cost of additional computation time.

\begin{table}[htbp]
\centering
\setlength{\tabcolsep}{6pt}
\resizebox{\textwidth}{!}{
\begin{tabular}{>{\centering\arraybackslash}p{1cm} >{\centering\arraybackslash}m{7cm} >{\centering\arraybackslash}m{7cm}}
\toprule
\textbf{Function} & \textbf{Precision $10^{-3}$} & \textbf{Precision $10^{-5}$} \\
\midrule
\addlinespace[5pt]  
$f_1$ &
\includegraphics[height=4.3cm]{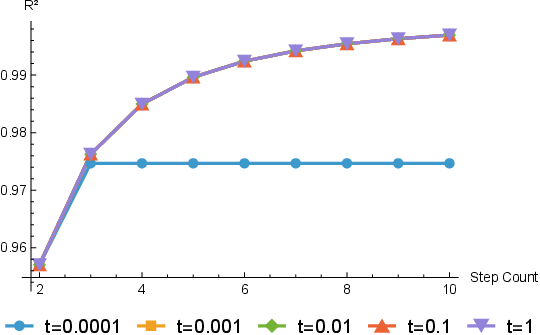} &
\includegraphics[height=4.3cm]{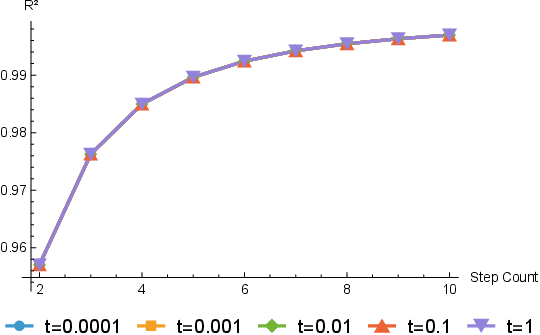} \\
\addlinespace[5pt]  
\midrule
\addlinespace[5pt]  
$f_3$ &
\includegraphics[height=4.3cm]{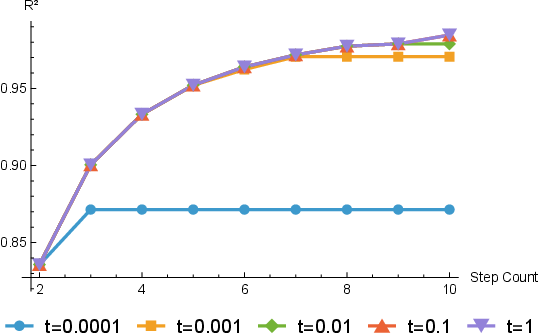} &
\includegraphics[height=4.3cm]{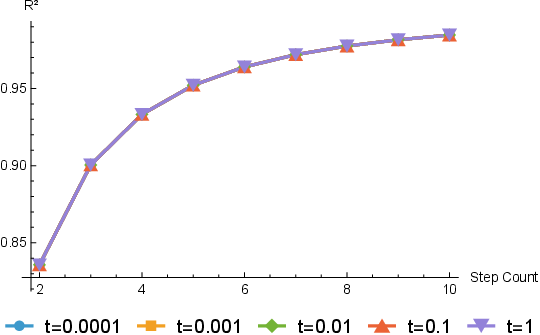} \\
\addlinespace[5pt]  
\midrule
\addlinespace[5pt]  
$f_6$ &
\includegraphics[height=4.3cm]{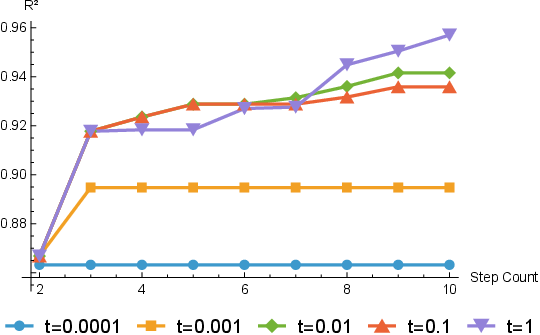} &
\includegraphics[height=4.3cm]{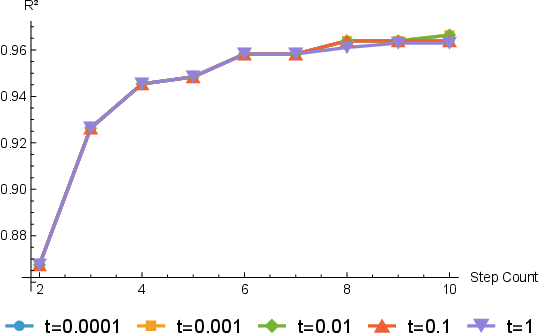} \\
\addlinespace[5pt]  
\midrule
\addlinespace[5pt]  
$f_7$ &
\includegraphics[height=4.3cm]{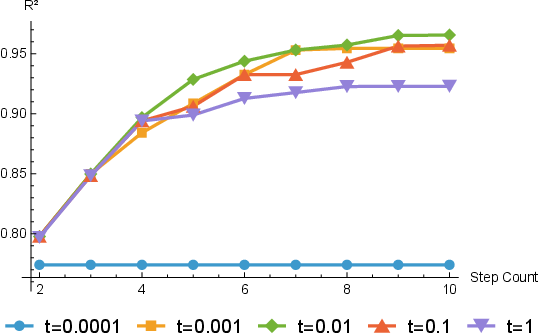} &
\includegraphics[height=4.3cm]{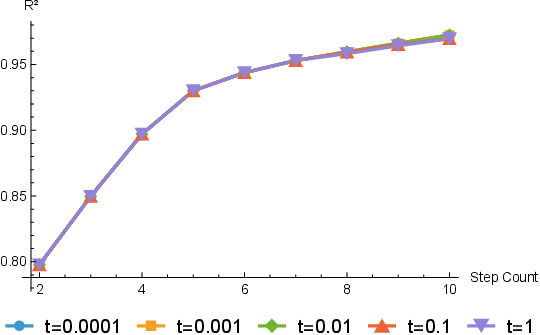} \\
\addlinespace[5pt]  
\bottomrule
\end{tabular}}
\caption{Comparison of $R^2$ tolerance variation under different precisions}
\label{tab:R2-tolerance-vs-precision}
\end{table}

\clearpage
\printbibliography

\appendix

\section{C++ code implementation}    \label{app:code}

\noindent\footnotesize{\tt{\# License:\ CC BY-NC 4.0:\ Attribution-NonCommercial 4.0 Int'l \href{https://creativecommons.org/licenses/by-nc/4.0/}{creativecommons.org/licenses/by-nc/4.0}}}
\lstdefinestyle{mystyle}{
    backgroundcolor=\color{backcolor},   
    commentstyle=\color{codegreen},
    keywordstyle=\color{blue},     numberstyle=\tiny\color{codegrey},
    stringstyle=\color{codepurple},
    basicstyle=\ttfamily\footnotesize,
    escapeinside={(*@}{@*)},
    breakatwhitespace=false,         
    breaklines=true,                 
    captionpos=b,                    
    keepspaces=true,                 
    numbers=left,                    
    numbersep=5pt,                  
    showspaces=false,                
    showstringspaces=false,
    showtabs=false,                  
    tabsize=2,
    emph={None},
    emphstyle={\color{magenta}}
}

\lstset{style=mystyle}

\begin{lstlisting}[language=C++,style=mystyle]
// EscalierAlgos.cpp
// (c) 2025 S. Bossu -- All rights reserved
(*@
\zigzag{15} @*)
// Global variables:
static double precision = 0, tolerance = 0;
static f_func f = nullptr;
static fbar_func fbar = nullptr;
(*@
\zigzag{15} @*)
// Local functions:
static double buildESS(double k0, double k1, double b, int nsteps);
static double getTail(double k0, double k1, double b, int nsteps, double*& reverse_fit);

/// <summary>
/// Find best escalier fit a < k1 < k2 < ... < b
/// </summary>
/// <param name="f">Function to be fitted</param>
/// <param name="fbar">Mean value function over an interval</param>
/// <param name="precision">Step position increment</param>
/// <param name="tolerance">Used for zero test of criticality</param>
/// <param name="maxsteps">Maximum number of steps</param>
/// <param name="a">Floor position</param>
/// <param name="b">Maximum step position</param>
/// <param name="reverse_fit">Output steps vector of size maxsteps in reverse order k[n] > k[n-1] > ... > k1</param>
/// <returns>ESS of solution</returns>
double EscalierFit(f_func f, fbar_func fbar, double precision, double tolerance, int maxsteps,
					double a, double b, double*& reverse_fit) {
(*@
\hspace{1.5em}\zigzag{14.5} @*)
	::precision = precision;
	::tolerance = tolerance;
	::f = f;
	::fbar = fbar;

	double k1, argmax = a; // step positions
	int i;
	double ess, maxess = 0;
	for (i = 0; (k1 = a + i * precision) < b; i++) {		// sweep search of first step position k1 with highest ess		
		ess = buildESS(a, k1, b, maxsteps);
		if (ess > maxess) {	// store best k1 in argmax and its ess in maxess
			maxess = ess;
			argmax = k1;
		}
	}

	reverse_fit = new double[maxsteps];
	ess = getTail(a, argmax, b, maxsteps, reverse_fit);	//
  (*@
\zigzag{14.5} @*)
	return maxess;
}

/// <summary>
/// Recursively compute ESS of step positions k0 < k1 < k2 < ... < b by sweep search of critical step k2
/// </summary>
/// <param name="k0">Floor position</param>
/// <param name="k1">First step position</param>
/// <param name="b">Maximum step position</param>
/// <param name="maxsteps">Maximum number of steps</param>
/// <returns>ESS value</returns>
double buildESS(double k0, double k1, double b, int maxsteps) {
	double maxess = (b - k1) * pow(fbar(k1, b), 2);		// base ESS if k2 = b (no further steps)
	if (maxsteps > 1) {
		double c, czone, argmin, k2; // c, czone: criticality of k2 (perfect = 0), argmin: optimal k2
		for (int i = 1; (k2 = k1 + i * precision) < b; i++) {	// sweep search of k2 between k1 and b
			c = abs(fbar(k0, k1) + fbar(k1, k2) - 2 * f(k1));	// criticality of k2
			if (c > EPSILON && c < tolerance) { // critical zone
				argmin = k2;
				do { // store best k2 in argmin
					i++; 
					if ((k2 = k1 + i * precision) >= b)
						break;
					czone = abs(fbar(k0, k1) + fbar(k1, k2) - 2 * f(k1));
					if (c > czone) {
						c = czone;
						argmin = k2;
					}
				} while (czone < tolerance);
				maxess = max(maxess, buildESS(k1, argmin, b, maxsteps - 1));  // Recursive call with one fewer step
			}
		}
	}

	return maxess + (k1 - k0) * pow(fbar(k0, k1), 2);
}

/// <summary>
/// Recursively compute the positions of optimal step positions k0 < k1 < k2 < ... < b by sweep search of critical step k2
/// </summary>
/// <param name="k0">Floor position</param>
/// <param name="k1">First step position</param>
/// <param name="b">Maximum step position</param>
/// <param name="maxsteps">Maximum number of steps</param>
/// <param name="reverse_fit">Output vector (kn, ..., k2, k1) of size n = maxsteps in reverse order</param>
/// <returns>ESS of optimal vector</returns>
double getTail(double k0, double k1, double b, int maxsteps, double*& reverse_fit) {
	double maxess = (b - k1) * pow(fbar(k1, b), 2);		// base ESS if k2 = b (no further steps)

	for (int i = 0; i < maxsteps; i++)	// initialize reverse_fit array
		reverse_fit[i] = b;
	if (maxsteps == 1)
		reverse_fit[0] = k1;
	else { // Case maxsteps > 1
		double c, czone, argmin, k2, ess; // c, czone: criticality of k2, argmin: optimal k2
		int k2count = 0;
		bool essIncreased = false;	// flag if any critical k2 identified has increased ESS
		for (int i = 1; (k2 = k1 + i * precision) < b; i++) {	// sweep search of k2 between k1 and b
			c = abs(fbar(k0, k1) + fbar(k1, k2) - 2 * f(k1));	// criticality of k2
			if (c > EPSILON && c < tolerance) { // critical zone
				argmin = k2;
				do { // store best critical k2 in argmin
					i++;
					if ((k2 = k1 + i * precision) >= b)
						break;
					czone = abs(fbar(k0, k1) + fbar(k1, k2) - 2 * f(k1));
					if (c > czone) {
						c = czone;
						argmin = k2;
					}
				} while (czone < tolerance);

				k2count++;
				if (k2count == 1) {	// first critical k2 found, use the same reverse_fit array
					ess = getTail(k1, argmin, b, maxsteps - 1, reverse_fit); // Recursive call with one fewer step
					if (maxess < ess) {
						essIncreased = true;
						maxess = ess;
						reverse_fit[maxsteps - 1] = k1;
					}
				}
				else {	// more than one critical k2 found, create new_reverse_fit array and copy results
					double* new_reverse_fit = new double[maxsteps];
					ess = getTail(k1, argmin, b, maxsteps - 1, new_reverse_fit); // Recursive call with one fewer step
					if (maxess < ess) {
						essIncreased = true;
						maxess = ess;
						reverse_fit[maxsteps - 1] = k1;
						for (int i = 0; i < maxsteps - 1; i++)
							reverse_fit[i] = new_reverse_fit[i];
					}
					delete[] new_reverse_fit;
				}
			}
		}
		if (!essIncreased) { // No critical k2 found with better ESS
			reverse_fit[maxsteps - 1] = k1;
			for (int i = 0; i < maxsteps - 1; i++)
				reverse_fit[i] = b;
		}
	}

	return maxess + (k1 - k0) * pow(fbar(k0, k1), 2);
}
(*@
\zigzag{15} @*)
\end{lstlisting}
    
\end{document}